\documentclass[reqno,12pt]{amsart}
\usepackage{amsmath, amssymb, amsthm, amsfonts} 
\usepackage[english]{babel}
\usepackage{bbm}
\usepackage{graphicx}
\usepackage{url}
\usepackage{soul}
\usepackage{epstopdf}
\usepackage{float}
\usepackage[ruled,vlined]{algorithm2e}
\usepackage[a4paper,bindingoffset=0.5cm,left=2cm,right=2cm,top=2.5cm,bottom=2cm,footskip=.8cm]{geometry}
\usepackage{rotating}
\usepackage{amsbsy,enumerate}
\usepackage{comment}
\usepackage{mathrsfs} 

\newcommand{\bs}{\boldsymbol}

\newcommand{\vb}{\vspace{3.2mm}}
\renewcommand{\hat}{\widehat}

\newcommand{\vertiii}[1]{{\left\vert\kern-0.25ex\left\vert\kern-0.25ex\left\vert #1 \right\vert\kern-0.25ex\right\vert\kern-0.25ex\right\vert}}
\allowdisplaybreaks
\usepackage{hyperref}

\newcommand{\jiesen}[1]{{\color{black}#1}}
\newcommand{\js}[1]{{\color{black}#1}} 
\newcommand{\michel}[1]{{\color{black}#1}}
\newcommand{\MM}[1]{{\color{black}#1}} 
\newcommand{\MRHM}[1]{{\color{black}#1}} 

\allowdisplaybreaks

\newtheorem{remark}{Remark}

\newtheorem{proposition}{Proposition}

\usepackage[utf8]{inputenc}
\usepackage{type1cm}         
\usepackage{multicol}        
\usepackage[bottom]{footmisc}

\usepackage{newtxtext}       %
\usepackage{newtxmath}       

\usepackage{tikz}
\usetikzlibrary{plotmarks}
\usetikzlibrary{decorations.pathreplacing}
\usepackage{pgfplots}
\pgfplotsset{width=10cm,compat=1.9}
\usepgfplotslibrary{external}
\usepackage{subcaption}
\usepackage{standalone}
\usetikzlibrary{automata,positioning}
\usetikzlibrary{decorations.pathreplacing,decorations.markings,shapes.geometric}
\usetikzlibrary{shapes,arrows}
\usetikzlibrary{backgrounds,calc,positioning}

\usepackage{ragged2e}
\usepackage{booktabs,tabularx}

\usepackage{pgfplots}
\pgfplotsset{compat=1.9}

\usetikzlibrary{chains,shapes.multipart}
\usetikzlibrary{shapes,calc}
\usetikzlibrary{automata,positioning}

\usepackage{tikz}
\usetikzlibrary{automata,positioning}
\usetikzlibrary{decorations.pathreplacing,decorations.markings,shapes.geometric}
\usetikzlibrary{shapes,arrows}
\usetikzlibrary{backgrounds,calc,positioning}

\begin{document}

	\title[Estimation in a dynamic Erd\H{o}s-R\'enyi random graph]{Estimation of on- and off-time distributions \\in a dynamic Erd\H{o}s-R\'enyi random graph}
\author{Michel Mandjes
{\tiny and} Jiesen Wang}
	
    \begin{abstract}
    In this paper we consider a dynamic Erd\H{o}s-R\'enyi graph in which edges, according to an alternating renewal process, change from present to absent and vice versa. The objective is to estimate the on- and off-time distributions while only observing the aggregate number of edges. This inverse problem is dealt with, in a parametric context, by setting up an estimator based on the method of moments. We provide conditions under which the estimator is asymptotically normal, and we point out how the corresponding covariance matrix can be identified. It is also demonstrated how to adapt the estimation procedure if alternative subgraph counts are observed, such as the number of wedges or triangles.

\vb


\vb

\noindent
{\sc Affiliations.} MM is with the Mathematical Institute, Leiden University, P.O. Box 9512,
2300 RA Leiden,
The Netherlands. He is also affiliated with Korteweg-de Vries Institute for Mathematics, University of Amsterdam, Amsterdam, The Netherlands; E{\sc urandom}, Eindhoven University of Technology, Eindhoven, The Netherlands; Amsterdam Business School, Faculty of Economics and Business, University of Amsterdam, Amsterdam, The Netherlands.


\noindent
JW is with the Korteweg-de Vries Institute for Mathematics, University of Amsterdam, Science Park 904, 1098 XH Amsterdam, The Netherlands.

\vb

\noindent
{\sc Acknowledgments.} 
The authors thank Camiel Koopmans (Mathematical Institute, Leiden University) and Liron Ravner (Department of Statistics, University of Haifa) for useful comments.
This research was supported by the European Union’s Horizon 2020 research and innovation programme under the Marie Sklodowska-Curie grant agreement no.\ 945045, and by the NWO Gravitation project NETWORKS under grant agreement no.\ 024.002.003. \includegraphics[height=1em]{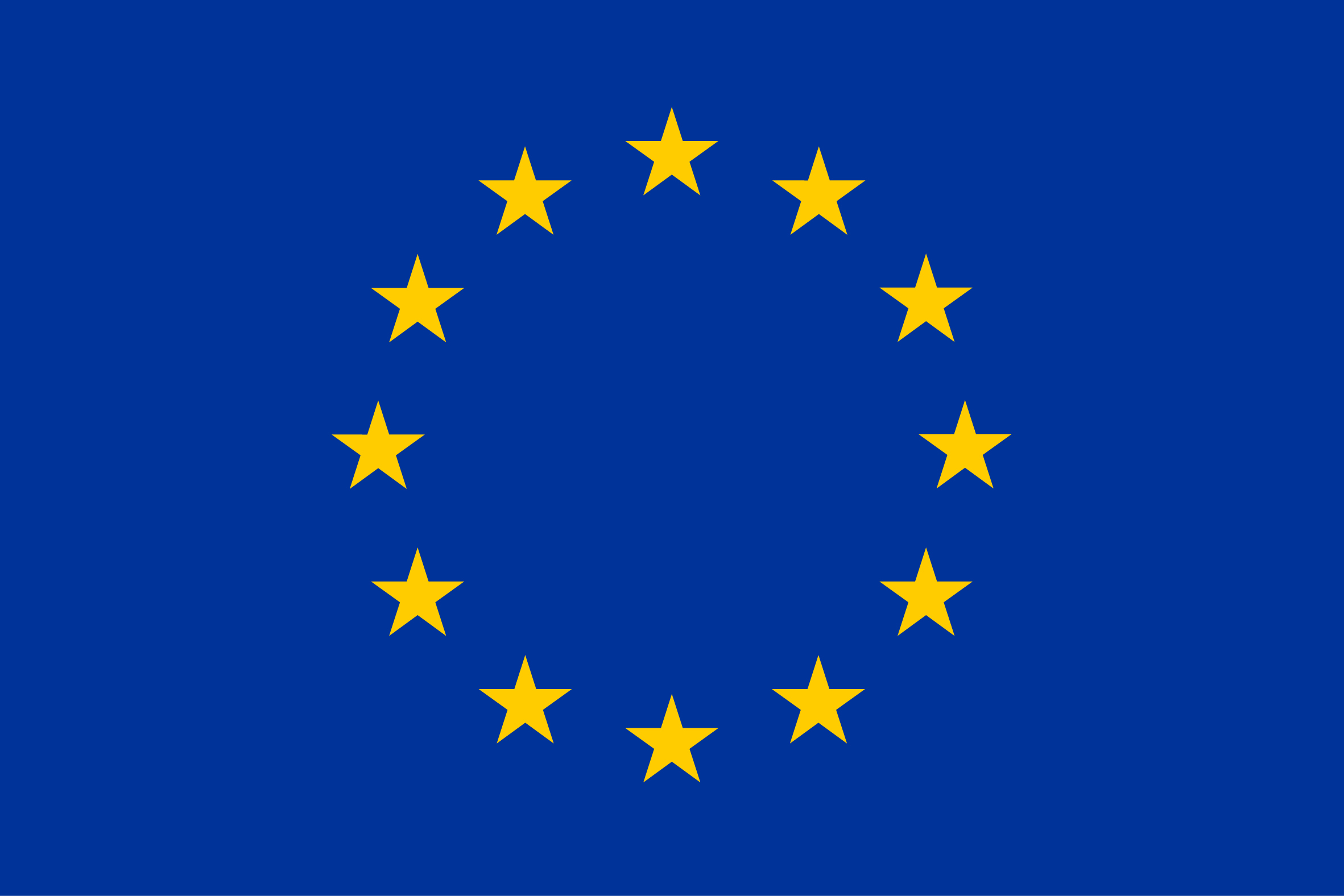} 
Date: {\it \today}.

\vb

\noindent
{\sc Email.} \url{m.r.h.mandjes@math.leidenuniv.nl}
and \url{j.wang2@uva.nl}.

	\end{abstract}

	\maketitle

 \newpage 
 
\section{Introduction}
{\it Model.} Consider a discrete-time model in which $n$ edges independently alternate between being present and being absent. 
The per-edge on-times form a sequence of iid (independent and identically distributed) random variables, where each on-time is distributed as a generic random variable $X\in{\mathbb N}$, while the off-times are iid and distributed as a generic random variable $Y\in{\mathbb N}$. 
It is in addition assumed that the sequence of on-times and the sequence of off-times are independent.
The resulting model can be seen as a dynamic version of the classical {\em Erd\H{o}s-R\'enyi random graph}: in case $n=\binom{N}{2}$ for some positive integer $N$ (to be interpreted as the number of vertices), each of the $n$ potential edges between vertex-pairs alternates between `on' and `off' according to the mechanism described above, independently of the other potential edges.
In equilibrium, the system behaves as a (static) Erd\H{o}s-R\'enyi random graph with $N$ vertices and `on-probability'
\[\varrho\equiv \varrho_{X,Y}:=\frac{{\mathbb E}\,X}{{\mathbb E}\,X+{\mathbb E}\,Y}.\]
where we throughout assume that ${\mathbb E}\,X$ and ${\mathbb E}\,Y$ are finite. 

\vb

{\it Objective.} In this paper we aim to devise a statistical inference technique by which we can estimate the distributions of the on-time $X$ and the off-time $Y$.
Crucially, in the context considered there is `partial information', in that the on- and off-times are {\it not} directly observed.
Indeed, letting ${\boldsymbol 1}_j(k):=1$ if edge $j\in\{1,\ldots,n\}$ is on at time $k\in{\mathbb N}$ and ${\boldsymbol 1}_j(k):=0$ otherwise,
the goal is to estimate the distribution of the on-times and off-times by only observing the {\it aggregate number of edges}
\[A_n(k):=\sum_{j=1}^n  {\boldsymbol 1}_j(k), \] 
at times $k=1,\ldots,K$ (with $K\in {\mathbb N}$).
Importantly, this means that we do not observe the evolution of individual edges, entailing that we do not have direct access to realizations of the distributions that we are interested in. 
In the most general setting, the objective is to estimate, for $k\in{\mathbb N}$,
\[F(k):={\mathbb P}(X\leqslant k),\:\:\:\:G(k):={\mathbb P}(Y\leqslant k),\]
or, equivalently,  the corresponding densities \[f_k:={\mathbb P}(X=k),\:\:\:\:g_k:={\mathbb P}(Y=k).\] 

Throughout this paper, we impose the natural assumption that the trace we are working with corresponds to a {\it stationary version} of the process. 
Concretely, this means that the system is in stationarity at time $1$: at that epoch, each edge is on with probability $\varrho$ and off otherwise, and the residual on-time $\bar X$ and the residual off-time $\bar Y$ are characterized by the densities, for $k\in{\mathbb N}$,
\[\bar f_k:= \frac{\sum_{\ell=k}^\infty f_\ell}{{\mathbb E}\,X},\:\:\:\bar g_k:= \frac{\sum_{\ell=k}^\infty g_\ell}{{\mathbb E}\,Y},\]
respectively. 

\vb

{\it Complications.}
In the case that $X$ and $Y$ are geometrically distributed, the process $\{A_n(k)\}_{k\in{\mathbb N}}$ is Markovian, which makes statistical inference relatively straightforward: one could first determine the transition probabilities of this Markov chain, and then estimate the parameters by applying maximum likelihood. 
For any other distribution of the on- and off-times, one cannot follow this convenient procedure, as $\{A_n(k)\}_{k\in{\mathbb N}}$ ceases to be Markovian. 
\michel{Observe that instances at which $A_n(k)$ equals 0 or $n$ are \js{not} regeneration points; for instance when $A_n(k)=0$, it matters how long each of the edges has been off already.}

\begin{remark}{\em 
    Technically, there are regeneration points that one could try to work with. Observe that, using the $k$ for which one has that $A_n(k) = 0$ and $A_n(k+1)=n$, one can estimate $f_1$; to this end, realize that in this scenario $A_n(k+2)$ is binomially distributed with parameters $n$ and $1-f_1$. 
    More generally, along the same lines one can argue that using the $k$ for which $A_n(k) = 0$ and $A_n(k+1)=\cdots=A_n(k+i-1)=n$, one can estimate $f_1+\cdots+f_{i-1}$. 
    In the same way, instances at which one goes from the full graph to the empty graph can be used to estimate the distribution of $Y$. 
    Note, however, that this approach is highly impractical: \michel{regeneration epochs are rare, in particular when the number of edges $n$ is large, rendering a substantial amount of information unused.} $\hfill\Diamond$}
\end{remark}

In this paper, we have cast our estimation problem in terms of edges of an evolving graph, but various alternative formulations can be chosen. 
In the context of a superposition of on-off traffic sources \cite[Chapter VIII]{WHI}, for instance, the estimation problem would amount to inferring the sources' on- and off-time distributions from observations of the aggregate rate. 

\vb

{\it Test distributions.} 
Our approach will be {\it parametric}: given parametric forms of the on- and off-time probability density functions, we will devise a method to estimate these underlying parameters. 
So as to cover a wide range of on- and off-times, the distributions we will be working with in this paper are the following.
In the first place, we consider geometric distributions. We write $Z\sim {\mathbb G}(p)$ for the geometric random variable with `success probability' $p$, for $p\in(0,1)$, when we mean
\[{\mathbb P}(Z\geqslant i) = (1-p)^{i-1},\:\:i=1,2,\ldots.\]
Second, we consider the more general class of Weibull distributions, covering tails that are lighter and heavier than geometric. We write $Z\sim{\mathbb W}(\lambda,\alpha)$, with $\alpha,\lambda>0$, to denote
\[{\mathbb P}(Z\geqslant i) = e^{-\lambda(i-1)^\alpha},\:\:i=1,2,\ldots,\]
so that ${\mathbb W}(-\log(1-p),1)\equiv {\mathbb G}(p)$. 
If $\alpha<1$, $Z$ is subexponential \cite[Section A5]{ASM2}; as can be verified directly, all moments of $Z$ are finite, but the moment generating function is infinite for any positive argument. 
In the third place, to model power-law tails, we consider a class of Pareto-type distributions (also sometimes referred to as the class of Lomax distributions). We write $Z\sim {\mathbb P}{\rm ar}(C,\alpha)$ to denote, for $C,\alpha>0$,
\[{\mathbb P}(Z\geqslant i) = \frac{C^\alpha}{(C+(i-1))^\alpha},\:\:i=1,2,\ldots;\]
to make sure ${\mathbb E}\,Z<\infty$ one should assume that $\alpha>1$. In the special case that $Z\sim {\mathbb P}{\rm ar}(1,\alpha)$, we have that ${\mathbb P}(Z\geqslant i) =i^{-\alpha}$. This class of distributions is subexponential for any $\alpha>0$, with ${\mathbb E}\, Z^k<\infty$ if and only if $k<\alpha.$

\vb

{\it Contributions.} The main contribution of this paper is a procedure by which one can estimate the parameters of the on- and off-times distributions. Conditions are provided under which the estimator, based on the method of moments, can be proven to be asymptotically normal.  
It is demonstrated how the corresponding covariance matrix can be identified. 
Interestingly, we point out how the estimator can be adapted to the situation in which, instead of the number of edges, other subgraph counts are observed, such as the number of triangles or wedges. 
We extensively discuss several variants of the setup considered in this paper, including an account of the model's continuous-time counterpart. \js{Appendix~A provides a detailed computation of the variance and covariance of the estimator when the on- and off-times follow a geometric distribution, while Appendix~B considers the setting in which the on- and off-times are Pareto distributed.}

\michel{At the methodological level, our paper combines known results from statistics (in particular the so called {\it delta method} \cite{VdV}) with explicit expressions for covariances of subgraph counts, so as to propose estimators and to prove properties thereof. A novelty of our work lies in opening up this class of inference problems for dynamic random graphs, with various opportunities for followup research.  }

\vb

{\it Literature.} The estimation problem discussed above lies at the interface of two active fields of research: in the first place there is the area of what could be called {\it inverse problems}, where one aims at estimating model primitives based on partial information, and in the second place there is the area of {\it dynamic random graphs}, where one aims at characterizing the probabilistic properties of graphs that randomly evolve over time. 
We proceed by a brief account of some of the relevant contributions and research directions in both fields.

\begin{itemize}
    \item[$\circ$]
Many types of inverse problems have been investigated. 
In a quintessential example, one has longitudinal data of the population size of a branching process, e.g.\ of Galton-Watson type, and one aims to estimate (features of) the offspring distribution; see e.g.\ the classical textbook \cite{GUT}. 
In another example, relating to an epidemiological context, one wishes to estimate the parameters of a stochastic SIR model by observing the number of infected individuals over time; see e.g.\ \cite{CAU} and the literature overview therein. 
In the economic literature, some papers focus on inverse problems in the context of financial data \cite{ACL, DUF}. 
Finally, we mention the stream of literature in which the objective is to infer the input of a storage system from periodic workload observations; see e.g.\ the survey \cite{ASA} and specific contributions \cite{HAN,RAV}.
\item[$\circ$] Where in the random graph literature the traditional focus is on the static case, recently attention has shifted to random graphs that stochastically evolve over time; see e.g.\ \cite{ATH,CRA}. The present paper is in the tradition of the model studied in \cite{BRA}, in which the edges' on- and off-times are exponentially distributed; the main result of \cite{BRA} is a sample-path large deviations principle, extending its counterpart \cite{CHA} for the static graph. 
\end{itemize}

{\it Organization.} 
This paper is organized as follows. 
In Section~\ref{PREL} we present helpful preliminaries that are used throughout the paper. 
Then Section~\ref{MOM} presents our approach based on the method of moments, including an assessment of the asymptotic normality of the resulting estimator. While most of the paper focuses on estimation based on observations of the number of edges, in Section~\ref{SUB} we point out how the estimation procedure extends to other subgraph counts. 
Numerical results, illustrating the performance of our estimator, are presented in Section~\ref{NUM}.
In Section~\ref{DISC} we provide a discussion and concluding remarks; this section in particular covers an account of the model's continuous-time counterpart. 

\section{Preliminaries} \label{PREL}
In this section we present a number of elementary properties of the model under study, which will be used extensively in later sections. 
A secondary goal of the section is to explain why a maximum-likelihood based estimation approach is, in our context, problematic.

\vb

When aiming to estimate the parameters of the on- and off-time distributions, a first possible approach could be based on maximum likelihood.
Such an approach requires the evaluation of the probability
\[\ell(n_1,\ldots,n_K):={\mathbb P}(A_n(1)=n_1,\ldots,A_n(K)=n_K),\]
with $n_k\in\{0,\ldots,n\}$ for $k=1,\ldots,K$. 
We can then maximize (the logarithm of) this $\ell(n_1,\ldots,n_K)$ over the parameter vector, thus yielding a maximum likelihood estimator. 

In the specific case the on- and off-times are both geometrically distributed, the process $\{A_n(k)\}_{k\in{\mathbb N}}$ is Markovian. 
More specifically, the probability $\ell(n_1,\ldots,n_K)$ allows an explicit evaluation.
Indeed, one can write
\[\ell(n_1,\ldots,n_K)=\michel{{\mathbb P}(A_n(1)=n_1)}\prod_{k=2}^{K} {\mathbb P}(A_n(k) = n_k\,|\, A_n(k-1)=n_{k-1}),\] 
where the transition probabilities ${\mathbb P}(A_n(k) = n_k\,|\, A_n(k-1)=n_{k-1})$ allow explicit evaluation.
Indeed, conditional on the value of $A_n(k)$, $A_n(k+1)$ can be written as the sum of two binomial quantities: with $X\sim {\mathbb G}(p)$ and $Y\sim {\mathbb G}(q)$,
\[A_n(k+1)\sim {\mathbb B}{\rm in}(A_n(k), 1-p)+{\mathbb B}{\rm in}(n-A_n(k),q), \]
with the two quantities on the right-hand side being independent. We conclude that in this case maximum likelihood can be directly applied in order to estimate the parameters $p$ and $q$. 

In all other cases, however, we do not have access to $\ell(n_1,\ldots,n_K)$, so that maximum likelihood cannot be (directly) applied. 
It is noted, though, that we can still set up a computational procedure that outputs the joint moment generating function of $({\bs 1}(1),\ldots,{\bs 1}(K))$, where ${\bs 1}(k)$ is the indicator function of an arbitrary edge being on at time $k$; this procedure is presented in Subsection \ref{SS:JMGF} below. 
As will be pointed out in Remark \ref{R2}, the joint moment generating function can be used to evaluate the joint distribution of $({\bs 1}(1),\ldots,{\bs 1}(K))$ (or any subset of these $K$ indicator functions), which is a quantity we will be working with in Subsection~\ref{SS:GEN}.

\subsection{Joint moment generating function}\label{SS:JMGF}
Before moving to the situation in which we start off in stationarity, we first focus on a single edge that just turned on \js{at time $0$}, and consider the generic random vector $({\boldsymbol 1}(1),\ldots,{\boldsymbol 1}(k))$ that records whether the edge is on at times $1,\ldots,K$. In this context, our objective is to evaluate the underlying moment generating function: in self-evident notation, we wish to calculate
\[M_+({\bs\theta}):={\mathbb E}_+ \exp\left(\sum_{k=1}^K \theta_k {\boldsymbol 1}(k)\right).\]
Likewise, for an edge that just turned off, our aim is to evaluate, again in self-evident notation,
\[M_-({\bs\theta}):={\mathbb E}_- \exp\left(\sum_{k=1}^K \theta_k {\boldsymbol 1}(k)\right).\]
\js{In the sequel, the following objects play a pivotal role: 
\[v_k({\bs\theta}):={\mathbb E}_{k,+} \exp\left(\sum_{i=k}^K \theta_i {\boldsymbol 1}(i)\right),\:\:\:\:w_k({\bs\theta}):={\mathbb E}_{k,-} \exp\left(\sum_{i=k}^K \theta_i {\boldsymbol 1}(i)\right),\]
where the subscript `$k,+$' (`$k,-$', resp.) means that the edge just went on (off, resp.) at time $k$.}
It is readily checked that these moment generating functions can be evaluated by solving the following system of linear equations:
\[
    v_k({\bs\theta}) = \sum_{\ell=1}^{K-k} f_\ell \,\exp\left(\sum_{i=k}^{k+\ell-1}\theta_i\right)\,w_{k+\ell}({\bs\theta})+
\sum_{\ell=K-k+1}^\infty f_\ell \,\exp\left(\sum_{i=k}^{K}\theta_i\right),
\]
and likewise
\[
    w_k = \sum_{\ell=1}^{K-k} g_\ell \,v_{k+\ell}({\bs\theta})+
\sum_{\ell=K-k+1}^\infty g_\ell .
\]
These are $2K$ linear equations in equally many unknowns. The system is uppertriangular, allowing for efficient recursive solution. Evidently,
$M_+({\bs\theta})= v_1({\bs\theta})$ and $M_-({\bs\theta})= w_1({\bs\theta})$.

Now we proceed by considering the situation that the on-off process starts in its stationary version.
Using the same reasoning as above, we thus have that
\[ M({\bs\theta}) =   \varrho\,M_+^{\rm (res)}({\bs\theta}) +(1-\varrho)
\log M_-^{\rm (res)}({\bs\theta}),
\]
where
\[M_+^{\rm (res)} ({\bs\theta})=\sum_{\ell=1}^{K-1} \bar f_\ell \,\exp\left(\sum_{i=1}^{\ell}\theta_i\right)\,w_{1+\ell}({\bs\theta})+
\sum_{\ell=K}^\infty \bar f_\ell \,\exp\left(\sum_{i=1}^{K}\theta_i\right),\]
and likewise
\[M_-^{\rm (res)} ({\bs\theta}) = \sum_{\ell=1}^{K-1} \bar g_\ell \,v_{1+\ell}({\bs\theta})+
\sum_{\ell=K}^\infty \bar g_\ell .\]

\begin{remark}\label{R2}{\em 
We finish our procedure to generate the joint moment generating function by pointing out that we can, at least in principle, explicitly convert $M({\bs\theta})$ into the joint distribution of the vector $({\bs 1}(1),\ldots,{\bs 1}(K))$.
To demonstrate this procedure, consider the case $K = 3$ (which we need in Subsection \ref{SS:GEN}). 
We wish to compute, in self-evident notation, the eight probabilities $p[0,0,0],p[0,0,1],\ldots,p[1,1,1]$.
To this end, we insert all eight elements from $\{-\infty,0\}^3$ into $M({\bs\theta})$. For instance, $M(-\infty,0,0)$ corresponds to the scenario that at the first time epoch the edge is off, while at the other two epochs the edge can be on or off: 
\[M(-\infty,0,0) = p[0,0,0]+p[0,1,0]+p[0,0,1]+p[0,1,1]. \] 
We thus obtain a non-singular system of eight linear equations in the equally many unknowns, which can be solved in the standard fashion. 
Observe that $M(0,0,0)=1$ corresponds to the requirement that the eight probabilities have to sum to $1$.} $\hfill\Diamond$
\end{remark}

\subsection{Saddlepoint approximation}
Above we concluded that, except in the case of geometric on- and off-times, we cannot explicitly compute the probability $\ell(n_1,\ldots,n_K)$.
One could still set up an {\it approximate} maximum likelihood estimation procedure though, following an approach advocated in \cite{DAV,GUN}.

To this end, we first introduce the {\it Legendre transform} (or the {\it convex conjugate}) of the joint cumulant generating function $\log M({\bs\theta})$: 
\[I({\bs n}):=\sup_{\bs\theta}\left({\bs\theta}^\top{\bs n}- n \log M({\bs\theta})\right), \]
where (in our context) ${\bs n}\in\{0,\ldots,n\}^K$.
Write ${\bs\theta}({\bs n})$ for the optimizing argument in the definition of $I({\bs n})$, for a given vector ${\bs n}$ (which is unique, as ${\bs\theta}\mapsto {\bs\theta}^\top{\bs n}-\log M({\bs\theta})$ is concave). 
Then the {\it saddle-point approximation} \cite{Butler, DAN}, particularly accurate when $n$ is large, states that 
\[\log \ell({\bs n})\equiv \log \ell(n_1,\ldots,n_K) \approx - J({\bs n}) :=-\frac{K}{2}\log(2\pi) - \frac{1}{2}\log s({\bs n})-\,I({\bs n}),\]
with
\[s({\bs n}) :=\left. {\rm det}\left(\frac{\partial^2}{\partial \theta_i\,\partial\theta_j}n\log M({\bs\theta})\right)\right|_{{\bs\theta} ={\bs\theta}({\bs n})}.\]
The idea of `approximate maximum likelihood' is to maximize $-J({\bs n})$ over the parameter vector. {As a simplification, one often just maximizes the `exponential part' $-I({\bs n})$ only.} A formal justification of this type of procedure has recently been provided in \cite{JES}. 

In the context of our estimation problem, the drawback of the approximate maximum likelihood approach lies in the computational effort. 
For a given value of the parameter vector, evaluation of $J({\bs n})$ or $I({\bs n})$ requires the maximization of ${\bs\theta}^\top{\bs n}-\log M({\bs\theta})$ over the $K$-dimensional vector ${\bs\theta}$, where it is noted that the computational effort required to compute $M({\bs\theta})$ is (as a consequence of the underlying uppertriangular structure) quadratic in $K$.
It turns out that, in particular for larger values of $K$, this computation can be rather time consuming, despite the fact that ${\bs\theta}\mapsto {\bs\theta}^\top{\bs n}-\log M({\bs\theta})$ is concave.
On top of this computation, there is the `outer loop': to numerically identify the (approximate) maximum likelihood estimator, these computations have to be done for a sequence of parameter vectors. 

Because of the computational challenges that we are facing, we propose to use an estimator based on the method of moments (to be discussed in the next section) which is of intrinsically lower numerical complexity. \js{A second key advantage of the use of the method of moments is that it allows us to work with other subgraph counts than just edges, as \MRHM{explained} in Section \ref{SUB}.}

\section{Method of moments}\label{MOM}
Recall that our objective is to estimate the parameters of the on- and off-time distributions from observations of the number of edges $A_n(1),\ldots, A_n(K).$
In this section we propose an estimator based on the method of moments. We first define this estimator, and then investigate its asymptotic normality.

\subsection{Estimator} \label{SS:EST} As a general principle, in the method of moments one needs to have as many `moment equations' as parameters. 
In our approach the moment equations are based on the quantities, for $\ell\in{\mathbb N}$,
\[S_{n,0}(k):=A_n(k)=\sum_{j = 1}^n{\boldsymbol 1}_j(k),\:\:\:\:S_{n,\ell}(k):=A_n(k)A_n(k+\ell)=\left(\sum_{j = 1}^n{\boldsymbol 1}_j(k)\right)\left(\sum_{j = 1}^n{\boldsymbol 1}_j(k+\ell)\right).\]
We proceed by providing explicit expressions for the expectations of these quantities, from which the moment equations can be identified. 
As we have assumed that the process is in equilibrium at time $1$, it remains in equilibrium at any time in $1,\ldots,K.$ 
We thus obtain, with $s_{n,\ell}:=\mathbb{E}\,S_{n,\ell}(k)$,
\begin{align*}
    s_{n,0}&= \sum_{j = 1}^n \mathbb{E}\left[{\boldsymbol 1}_j(k)\right] = n\,  {\mathbb P}({\boldsymbol 1}_i(k)=1)=n\varrho.
    \end{align*}
Along the same lines, using the independence between the per-edge processes, 
    \begin{align*}
    s_{n,1}&= \mathbb{E} \left[\sum_{j = 1}^n {\boldsymbol 1}_j(k) {\boldsymbol 1}_j(k+1)\right] + \mathbb{E}\left[\sum_{j = 1}^n\sum_{i \ne j}{\boldsymbol 1}_j(k){\boldsymbol 1}_i(k+1) \right]\\
    &= n\,
    {\mathbb P}({\boldsymbol 1}_i(k)=1)\,{\mathbb P}({\boldsymbol 1}_i(k+1)=1\,|\,{\boldsymbol 1}_i(k)=1) + (n^2-n)\big( {\mathbb P}({\boldsymbol 1}_i(k)=1)\big)^2\\
    &= n \varrho\, (1-\bar f_1)  + (n^2-n) \, \varrho^2.
\end{align*}
Similar to the computation of $s_{n,1}$, we find
\begin{align}s_{n,2} &=n \varrho\, \big((1-\bar f_1-\bar f_2) +\bar f_1g_1\big) + (n^2-n) \, \varrho^2,\label{sn2}\\
s_{n,3}&= n \varrho\, \big((1-\bar f_1-\bar f_2-\bar f_3) +\bar f_1g_2 + \bar f_1 g_1 (1- f_1)+\bar f_2 g_1\big) + (n^2-n) \, \varrho^2.\notag
\end{align}
To understand the expression for $s_{n,2}$, it should be borne in mind that there are two scenarios in which the edge is on at time $k$ and $k+2$, namely (in self-evident notation) +++ and +$-$+. Regarding $s_{n,3}$, a similar reasoning applies but now there are four scenarios: ++++, +$-$$-$+, +$-$++, ++$-$+. Moment equations for $s_{n,\ell}$ with $\ell= 4,5,\ldots$ can be identified in the same fashion (where the enumeration of all possible scenarios becomes increasingly tedious). 

Let $\hat\mu_{n,K}(\ell)$ be an estimator of $s_{n,\ell}$, where evident choices are (with $\ell\in{\mathbb N}$)
\[\hat\mu_{n,K}(0):=\frac{1}{K}\sum_{k=1}^K A_n(k),\:\:\:
\hat\mu_{n,K}(\ell):=\frac{1}{K-\ell}\sum_{k=1}^{K-\ell} A_n(k)A_n(k+\ell).
\]
In a model with $L\in{\mathbb N}$ parameters, our method of moments estimator of these parameters is a parameter vector solving the $L$ moment equations
\[s_{n,\ell}=\hat\mu_{n,K}(\ell),\:\:\:\:\ell=0,\ldots, L-1. \]

We proceed by considering a number of special cases, which will return in the numerical experiments presented in Section \ref{NUM}. \michel{In general, the `on-probability' $\varrho$ can be estimated by $\hat{\varrho}_{n,K} = \hat{\mu}_{n,K}(0)/n$.} 

\begin{itemize}
    \item[$\circ$]
In the case that $X\sim{\mathbb G}(p)$ and $Y\sim{\mathbb G}(q)$, there are two parameters that need to be estimated. The method of moments therefore involves $s_{n,0}$ and $s_{n,1}$. 
The associated moment equations reduce to
\begin{align}
    s_{n,0}&= n\varrho = n\frac{q}{p+q},\:\:\:\:\:
    s_{n,1} =n\varrho\,(1-p) +(n^2-n)\varrho^2.\label{ME_GG}
\end{align}
We conclude that
\[\hat\varrho_{n,K}=\frac{\hat\mu_{n,K}(0)}{n},\:\:\:\:\:\hat p_{n,K}=1-\frac{\hat\mu_{n,K}(1)-(n^2-n)\hat\varrho_{n,K}^2}{n\hat\varrho_{n,K}},\]
so that our estimators become
\begin{align}\notag
 \hat p_{n,K}&=\frac{\hat\mu_{n,K}(0)-\hat\mu_{n,K}(1)+(1-1/n)(\hat\mu_{n,K}(0))^2}{\hat\mu_{n,K}(0)},\\
\hat q_{n,K}&=\frac{\hat\mu_{n,K}(0)-\hat\mu_{n,K}(1)+(1-1/n)(\hat\mu_{n,K}(0))^2}{n-\hat\mu_{n,K}(0)}.\label{qhat}
\end{align}

 \item[$\circ$]
Also in the case that $X\sim {\mathbb P}{\rm ar}(1,\alpha)$ and $Y\sim {\mathbb P}{\rm ar}(1,\beta)$, two moment equations are needed. 
Observe that ${\mathbb E}\,X=\zeta(\alpha)$ and ${\mathbb E}\,Y=\zeta(\beta)$, where $\zeta(\alpha): = \sum_{i = 1}^{\infty}1/i^{\alpha}$ (for $\alpha>1$) is the {\it Riemann zeta function}.
The moment equations thus become
\begin{align*}
    s_{n,0} &= n  \varrho = n \, \frac{\zeta(\alpha)}{\zeta(\alpha)+\zeta(\beta)} ,\:\:\:\:
    s_{n,1} = n\varrho\left(1-\frac{1}{\zeta(\alpha)} \right)+ (n^2-n) \varrho^2.
\end{align*} 
It follows, after some elementary \michel{calculations} that
\[\hat\varrho_{n,K}=\frac{\hat\mu_{n,K}(0)}{n},\:\:\:\:\:\hat \alpha_{n,K}={\zeta^{-1}}\left(\frac{\hat\mu_{n,K}(0)}{\hat\mu_{n,K}(0)-\hat\mu_{n,K}(1)+(1-1/n)(\hat\mu_{n,K}(0))^2}\right),\]
\michel{and similarly}
\begin{align*}
  \hat\beta_{n,K}&={\zeta^{-1}}\left(\frac{n-\hat\mu_{n,K}(0)}{\hat\mu_{n,K}(0)-\hat\mu_{n,K}(1)+(1-1/n)(\hat\mu_{n,K}(0))^2}\right).
\end{align*}

\item [$\circ$]

In the case that $X\sim 
{\mathbb W}(1,\alpha)$ and $Y\sim {\mathbb G}{\rm eo}(q)$, entailing that there are two parameters, two moment equations are needed. These are
\begin{align*}
    s_{n,0} &= n  \varrho = n \, \frac{\chi(\alpha)}{\chi(\alpha)+1/q} ,\:\:\:\:
    s_{n,1} = n\varrho\left(1-\frac{1}{\chi(\alpha)} \right)+ (n^2-n) \varrho^2\,,
\end{align*}
where $\chi(\alpha) := \sum_{i=1}^{\infty} e^{-(i-1)^\alpha}$. It follows that
$\hat q_{n,K}$ is given by \eqref{qhat} and
\begin{align*}
    & \hat \alpha_{n,K} = \chi^{-1} \left(\frac{\hat\mu_{n,K}(0)}{\hat\mu_{n,K}(0)-\hat\mu_{n,K}(1)+(1-1/n)(\hat\mu_{n,K}(0))^2}\right)
\end{align*}

 \item[$\circ$]
In the case that $X\sim {\mathbb P}{\rm ar}(C,\alpha)$ and $Y\sim {\mathbb G}{\rm eo}(q)$ there are three parameters, so that three moment equations are needed:
\begin{align*}
    s_{n,0} &= n  \varrho = n \, \frac{\zeta(C,\alpha)}{\zeta(C,\alpha)+1/q} ,\:\:\:\:
    s_{n,1} = n\varrho\left(1-\frac{1}{\zeta(C,\alpha)} \right)+ (n^2-n) \varrho^2\,, \\
    s_{n,2} &= n\varrho\left(1-\frac{1}{\zeta(C,\alpha)}-\frac{(C/(C+1))^{\alpha}}{\zeta(C,\alpha)} + \frac{q}{\zeta(C,\alpha)}\right) + (n^2-n) \varrho^2 \\
    &= s_{n,1} + n\varrho\left(\frac{q}{\zeta(C,\alpha)}-\frac{(C/(C+1))^{\alpha}}{\zeta(C,\alpha)} \right) \,,\:\:\:\:\:\:\:
    \zeta(C,\alpha):=\sum_{i=1}^\infty \frac{C^\alpha}{(C+(i-1))^\alpha}.
\end{align*}
As it turns out, we again have that 
$\hat q_{n,K}$ is given by \eqref{qhat}.
The values of $\hat C_{n,K}$ and $\hat \alpha_{n,K}$ are obtained by solving two equations in two unknowns:
\begin{align*}
\zeta(C,\alpha)&=\frac{\hat\mu_{n,K}(0)}{\hat\mu_{n,K}(0)-\hat\mu_{n,K}(1)+(1-1/n)(\hat\mu_{n,K}(0))^2} \\
    \left(\frac{C}{C+1}\right)^\alpha &= \hat q_{n,K}-
    \frac{\hat\mu_{n,K}(2)-\hat\mu_{n,K}(1)}{\hat\mu_{n,K}(0)-\hat\mu_{n,K}(1)+(1-1/n)(\hat\mu_{n,K}(0))^2}.
\end{align*}
 
\end{itemize}

\begin{remark}\label{RR1}\michel{\em 
     It is clear that $s_{n,0}$ depends on the distributions of $X$ and $Y$ through  ${\mathbb E}\,X$ and ${\mathbb E}\,Y$ only. Interestingly, this holds for $s_{n,1}$ as well: noting that $\bar{f}_1 = 1/{{\mathbb E}\,X}$, we see that also $s_{n,1}$ is a function of ${\mathbb E}\,X$ and ${\mathbb E}\,Y$. This means that in instances in which the on- and off-times both have one parameter, one can equivalently estimate ${\mathbb E}\,X$ and ${\mathbb E}\,Y$ and then find the parameter values that are in line with these values. 
    
    In this paper we primarily focus on on- and off-times with relatively few parameters. The number of moment equations should clearly match the number of parameters to be estimated. This means that a complication lies in the fact that for larger values of $k$ the covariance between $A_n(1)$ and $A_n(k)$ becomes small, so that $s_{n,k}$ `provides little information' (i.e.,  $s_{n,k}$ will become increasingly hard to distinguish from $n^2\varrho^2$ as $k$ grows). Consequently, in practical terms, there are evident limits to the number of parameters our method can handle.
    
    In the cases in which both $X$ and $Y$ are characterized through a single parameter, we work with moment equations based on $s_{n,0}$ and $s_{n,1}$. In principle, however, we could have as well set up an estimator based on, say, $s_{n,0}$ and $s_{n,2}$. Above we already pointed out a drawback of that choice: for larger $k$, $s_{n,k}$ becomes less informative. In addition, as indicated by the involved expression for $s_{n,2}$ given in \eqref{sn2}, the moment equations are less convenient to work with. 
    $\hfill\Diamond$}
\end{remark}

\subsection{Asymptotic normality for geometric on- and off-times} 
We proceed by establishing asymptotic normality for the case that $X\sim{\mathbb G}(p)$ and $Y\sim{\mathbb G}(q)$, applying the {\it delta method}; see e.g.\  \cite[Section 3.9]{VdV}. The main idea is that under the proviso that $(\hat \mu_{n,K}(0), \hat \mu_{n,K}(1))$ is asymptotically normal (as $K\to\infty$, that is), then this property is inherited by $(\hat p_{n,K}, \hat q_{n,K})$, by virtue of the delta method.
\js{The starting point is that we assume the vector $(\hat \mu_{n,K}(0), \hat \mu_{n,K}(1))$ satisfies the following convergence in distribution: as $K\to\infty$,}
\begin{equation}\label{asnmu}
  \sqrt{K}\left(\begin{array}{c}\hat \mu_{n,K}(0)-s_{n,0}\\\hat \mu_{n,K}(1)-s_{n,1}\end{array}
\right) \js{\overset{\rm d}{\to}} \, \MRHM{\bs Z}:=\left(\begin{array}{c}Z_{0}\\ Z_{1}\end{array}\right),
\end{equation}
\MRHM{with ${\boldsymbol Z}$ being a bivariate zero-mean normal vector};
\js{we justify this claim later in this subsection. At this point one can immediately apply the results from \cite[Section 3.9]{VdV}, via which one can establish the joint central limit theorem for  $(\hat p_{n,K}, \hat q_{n,K})$ by evaluating a matrix of partial derivatives, corresponding to the mapping between $(\hat \mu_{n,K}(0), \hat \mu_{n,K}(1))$ and $(\hat p_{n,K}, \hat q_{n,K})$. 
The claim of Proposition~\ref{P1} below is now a direct consequence of the general result presented in \cite[Thm.\ 3.9.4]{VdV}.
However, we choose to also include an intuitive, informal calculation here, as it reveals how properties of $(\hat \mu_{n,K}(0), \hat \mu_{n,K}(1))$ carry over to properties of $(\hat p_{n,K}, \hat q_{n,K})$.}

The informal reasoning is as follows. 
\begin{itemize}
    \item[$\circ$]
    \MRHM{The distributional limit \eqref{asnmu} can be written equivalently as}
\begin{equation}\MRHM{\label{R2sEquation(4.5)}}
  \left(\begin{array}{c}\hat \mu_{n,K}(0)\\\hat \mu_{n,K}(1)\end{array}
\right) \js{\overset{\rm d}{=}} \left(\begin{array}{c}s_{n,0}\\ s_{n,1}\end{array}\right) + \frac{1}{\sqrt{K}} 
\left(\begin{array}{c}Z_{0}\\ Z_{1}\end{array}\right)\MM{\,+\,o_{\mathbb P}\left(\frac{1}{\sqrt{K}}\right)},
\end{equation}
as $K\to\infty$.
\MRHM{The scaling in \eqref{asnmu}--\ref{R2sEquation(4.5)}, in combination with the definitions of $\hat \mu_{n,K}(0)$, $\hat \mu_{n,K}(1)$, and $s_{n,0}$, $s_{n_1}$ suggests that the random vector ${\bs Z}$, with mean zero, should have a covariance matrix that is given by the asymptotic rescaled variances and covariances}
 \begin{align}v_0 &:=\lim_{K\to\infty}\frac{1}{K}{\mathbb V}{\rm ar}\left[\sum_{k=1}^K A_n(k)\right],\:\: \notag
    v_1 :=\lim_{K\to\infty}\frac{1}{K-1}{\mathbb V}{\rm ar}\left[\sum_{k=1}^{K-1} A_n(k)A_n(k+1)\right],\:\:\\c_{01}&:=\lim_{K\to\infty}\frac{1}{K-1}{\mathbb C}{\rm ov}\left[\sum_{k=1}^{K} A_n(k), \sum_{k=1}^{K} A_n(k)A_n(k+1)\right].\label{covvar}\end{align}
\MRHM{Proposition \ref{P1} below confirms this intuition and} computes these $v_0$, $v_1$, and $c_{01}$ for the case of $X\sim{\mathbb G}(p)$ and $Y\sim{\mathbb G}(q)$, relying on the definitions in the previous display; in Subsection \ref{SS:GEN} we point out how the three quantities in \eqref{covvar} can be calculated in general, and we in particular show that they are finite in a
heavy-tailed example.
    \item[$\circ$]
\MM{As a second step, based} on Taylor approximations of $(\hat \mu_{n,K}(0), \hat \mu_{n,K}(1))$ around $(s_{n,0},s_{n,1})$, we write
\[\hat p_{n,K}= \frac{s_{n,0}+Z_0/\sqrt{K}-(s_{n,1}+Z_1/\sqrt{K})+(1-1/n)(s_{n,0}+Z_0/\sqrt{K})^2}{s_{n,0}+Z_0/\sqrt{K}}\MM{\,+\,o_{\mathbb P}\left(\frac{1}{\sqrt{K}}\right)},\]
and an analogous expression for $\hat q_{n,K}$.
From the moment equations \eqref{ME_GG} it follows immediately that
\[\frac{s_{n,0}-s_{n,1}+(1-1/n)s_{n,0}^2}{s_{n,0}}=p.\]
Upon combining this identity with the familiar linear approximation $(1+x)^{-1} = 1 - x +O(x^2)$ 
as $x\downarrow 0$, we obtain that
\begin{align*}\hat p_{n,K} &= p  +\frac{Z_0}{\sqrt{K}}\left(\frac{s_{n,1}}{s_{n,0}^2}+ 1-\frac{1}{n}\right)+ \frac{Z_1}{\sqrt{K}}\left(-\frac{1}{s_{n,0}}\right)\MM{\,+\,o_{\mathbb P}\left(\frac{1}{\sqrt{K}}\right)}.\end{align*}
\item[$\circ$]We thus obtain that $\hat p_{n,K}$ is asymptotically normally distributed around $p$, and a similar relation for $\hat q_{n,K}$: for vectors ${\bs\gamma}, {\bs\delta} \in {\mathbb R}^2$,
\begin{align*}\hat p_{n,K} &= p + \frac{1}{\sqrt{K}}\left(\gamma_{0}Z_0+\gamma_{1}Z_1\right)\MM{\,+\,o_{\mathbb P}\left(\frac{1}{\sqrt{K}}\right)},\\
\hat q_{n,K} &= q + \frac{1}{\sqrt{K}}\left(\delta_{0}Z_0+\delta_{1}Z_1\right)\MM{\,+\,o_{\mathbb P}\left(\frac{1}{\sqrt{K}}\right)}.\end{align*}
It takes a minor computation to verify that ${\bs \delta} = (q/p){\bs \gamma}$, where the entries of ${\bs\gamma}$ equal
\begin{equation}
    \label{defgg}
\gamma_0=\frac{s_{n,1}}{s_{n,0}^2}+ 1-\frac{1}{n} =
2+\frac{p-q-(p+q)p}{nq},\:\:\:\:\:\gamma_1=-\frac{p+q}{nq}.\end{equation}
In this way we have identified the parameters pertaining to the limiting distribution of $(\hat p_{n,K}, \hat q_{n,K})$. 
\end{itemize}
\MM{The following result presents the asymptotic normality of our estimator, in the regime $K\to\infty.$}

\begin{proposition}\label{P1}
Consider the case that $X\sim{\mathbb G}(p)$ and $Y\sim{\mathbb G}(q)$. \MRHM{Define, with $\gamma_0$ and $\gamma_1$ as given in \eqref{defgg},
\begin{equation}
    \sigma^2:=\gamma_0^2v_0+2\gamma_0\gamma_1c_{01}+\gamma_1^2v_1,\quad \tau^2:=\sigma^2\,(q/p)^2,\quad \rho:=\sigma^2 (q/p).\label{sigdef}
\end{equation}}
    \MRHM{Then the following claims apply. \begin{itemize}\item[(i)] The distributional convergence \eqref{asnmu} is valid: as $K\to\infty$,
    \[\sqrt{K}\left(\begin{array}{c}\hat \mu_{n,K}(0)-s_{n,0}\\\hat \mu_{n,K}(1)-s_{n,1}\end{array}
\right) \overset{\rm d}{\to}
{\mathscr N}\left(\Big(\begin{array}{c}0\\0\end{array}\Big), \Big(\begin{array}{cc}v_0&c_{01}\\
c_{01}&v_1
\end{array}\Big)\right),\]
with $v_0$, $v_1$, and $c_{01}$ as defined in \eqref{covvar}.
 \item[(ii)] As $K\to\infty$,
\[\sqrt{K}\left(\begin{array}{c}
\hat p_{n,K} -p\\
\hat q_{n,K} -q\end{array}
\right)\to 
 {\mathscr N}\left(\Big(\begin{array}{c}0\\0\end{array}\Big), \Big(\begin{array}{cc}\sigma^2&\rho\\
\rho&\tau^2
\end{array}\Big)\right)
,\]
{with $\sigma^2$, $\tau^2$, and $\rho$ given in \eqref{sigdef}}.
   \item[(iii)] In addition,
    \begin{equation} \label{eq:nu0}
        v_0= n\,\varrho(1-\varrho)\frac{1+\lambda}{1-\lambda}=n\,\frac{pq(2-p-q)}{(p+q)^3},
    \end{equation}
    and $v_1$ is given by \eqref{v1} and $c_{01}$ by \eqref{c01}.\end{itemize}}
\end{proposition}

\js{{\it Proof.} \MRHM{Noting that above we have established part (ii)}  under the assumption that part (i) applies,
we \MRHM{now verify} the validity of \eqref{asnmu}.
In the case of $X \sim {\mathbb G}(p)$, $Y \sim {\mathbb G}(q)$, we have that $\{A_n(k),A_n(k+1))\}$ is an irreducible Markov chain, on the finite state-space $\{0,\ldots,n\}\times\{0,\ldots,n\}$, that we have assumed to be in stationarity. As a consequence, a standard central limit theorem for irreducible finite-state Markov chains in stationarity, states that 
 \[
\sqrt{K}\left(\frac{1}{K} \sum_{k=1}^K f(A_n(k), A_n(k+1)) - {\mathbb E} f(A_n(1), A_n(2))\right)
 \]
 converges in distribution, as $K\to\infty$, to a centered Normal random variable with a variance that depends on the function $f$ chosen, namely
 \[ \sigma^2_f:={\mathbb V}{\rm ar} f(A_n(1), A_n(2)) + 2\sum_{k=1}^\infty {\mathbb C}{\rm ov}\big(f(A_n(1), A_n(2)),f(A_n(k+1), A_n(k+2))\big).\]
 \MRHM{Observe that $v_0$ and $v_1$ from \eqref{covvar} correspond to $\sigma^2_f$ for the cases $f(x,y)=x$ and $f(x,y)=xy$, respectively.}
 Supposing that we pick $f(x,y):=\alpha x + \beta xy$ for arbitrary scalars $\alpha,\beta$,
 \MRHM{we find that $\sigma^2_f = \alpha^2\,v_0+\beta^2\, v_1+2\alpha\beta\, c_{01}$}, and that
 \[\sqrt{K} \Big(\alpha\hat \mu_{n,K}(0)+\beta\hat \mu_{n,K}(1)-\alpha s_{n,0}-\beta s_{n,1}\Big) \]
 converges to a centered Normal random variable, so that with the \MM{Cram\'er-Wold device \cite[Thm.\ 29.4]{BILL}} it follows that  
  \[\sqrt{K} \Big(\hat \mu_{n,K}(0)-  s_{n,0},\hat \mu_{n,K}(1)- s_{n,1}\Big) \]
  converges to a centered bivariate Normal random vector. \MRHM{It is noted that} the finiteness of $v_0$, $v_1$ and $c_{01}$ is guaranteed by part (iii) below. }

\MRHM{We now turn to part (iii), demonstrating \MM{how the quantities introduced in \eqref{covvar}}  can be computed.}
As we are in the setting of $X\sim{\mathbb G}(p)$ and $Y\sim{\mathbb G}(q)$, 
we can exploit the Markovian nature of the process $\{A_n(k)\}_{k\in{\mathbb N}}$. \MM{Define $\lambda:=1-p-q.$}
We focus here on the computation of $v_0$.
The quantities $v_1$ and $c_{01}$ can be dealt with in a similar way, but require considerably more intricate calculations; these are provided in Appendix~\ref{AppA}. 
\MRHM{We are to evaluate}
\begin{align*}
        \MRHM{v_0} &= \MRHM{{\mathbb V}\mbox{ar}\,A_n(1) + 2 \sum_{k=2}^\infty {\mathbb C}\mbox{ov}\left(A_n(1),A_n(k)\right).}
    \end{align*}
    Now notice that, conditional on $A_n(1)$, with the two binomial random variables appearing in the right-hand side being independent,
 \[A_n(k) \sim {\mathbb B}{\rm in}(A_n(1), r_{k-1}) + {\mathbb B}{\rm in}(n-A_n(1),s_{k-1}),\]
    where, as follows by an elementary computation,
    \begin{align*}r_k&:={\mathbb P}({\boldsymbol 1}_j(k)=1\,|\,{\boldsymbol 1}_j(1)=1)=\varrho +(1-\varrho)\,\lambda^{k-1},\\
    s_k&:={\mathbb P}({\boldsymbol 1}_j(k)=1\,|\,{\boldsymbol 1}_j(1)=0)=\varrho -\varrho\,\lambda^{k-1}.
    \end{align*}
    We thus find that
    \begin{align*}
       {\mathbb E}\big[A_n(1)A_n(k)\big] &= {\mathbb E}\big[A_n(1)\big(A_n(1)\,r_{k-1} + (n-A_n(1))\,s_{k-1}\big)\big] = 
    \lambda^{k-1}\,n\varrho(1-\varrho)+n^2\varrho^2,
    \end{align*}
    so that ${\mathbb C}\mbox{ov}\left(A_n(1),A_n(k)\right)=\lambda^{k-1}\,n\varrho(1-\varrho).$
    \MRHM{Upon combining the above, and in addition using that ${\mathbb V}\mbox{ar}\,A_n(1)=n\,\varrho(1-\varrho)$, we obtain that
    \[v_0 = n\,\varrho(1-\varrho) + 2n\,\varrho(1-\varrho)\sum_{k=2}^{\infty}\,\lambda^{k-1},\]
    which yields the stated expression.}
   \MM{As mentioned, the computations for $v_1$ and $c_{01}$ are presented in Appendix \ref{AppA}.}
\hfill$\Box$

\js{
\begin{remark}\label{RPQ}\em In case $\lambda = 0$, i.e., $p+q=1$, it is readily verified that that $\{A_n(k)\}_{k\in{\mathbb N}}$ is a sequence of {\it independent} ${\mathbb B}{\rm in}(n, \varrho)$ random variables. This could suggest that one can only estimate $\varrho$, and not the individual values of $p$ and $q$, but this intuition is wrong: the proposed estimators $(\hat p_{n,K},\hat q_{n,K})$ {\it are} capable of recuperating the true values $p$ and $q$ in an asymptotically Normal fashion. The crucial realization is that from the data, as $K$ grows large, it is increasingly picked up that we are dealing with an instance in which $s_{n,1}$ equals $s_{n,0}^2$. This means that, by \eqref{ME_GG},
\[n^2\varrho^2=n\varrho\,(1-p) +(n^2-n)\varrho^2,\]
which can be checked to reduce to $p+q=1$. Being able to estimate $\varrho$, together with knowing that $p+q=1$, makes the individual parameters $p$ and $q$ identifiable. Informally speaking, the fact that we learn that there is no correlation between subsequent observations, does provide helpful information that can be used in the estimation. 
\hfill$\Diamond$
\end{remark}}

\subsection{General recipe for computing (co-)variances} \label{SS:GEN}
Whereas in the above analysis for the case of geometric on-and off-times we intensively relied on the memoryless property of the geometric distribution, rendering $\{A_n(k)\}_{k\in{\mathbb N}}$ Markov, our next goal is to set up an approach to compute $v_0$, $v_1$, and $c_{01}$ for {\it any} pair of on- and off-time distributions. 
In the analysis of the case with geometric on- and off-times (see in particular the computations for $v_1$ in Appendix \ref{AppA}), we have observed that the most challenging component in the computation is ${\mathbb E}[A_n(1)A_n(2)A_n(k)A_n(k+1)]$, with the system being in stationarity at time $1$. For that reason, 
we demonstrate how this object can be evaluated for general on- and off-times, leaving it to the reader to adapt this procedure to make it applicable to all objects appearing in $v_0$, $v_1$, and $c_{01}$. 

The analysis starts by defining the following probabilities: for $i_2,i_k,i_{k+1}\in\{0,1\}$,
\begin{align*}
p^-_{i_2,i_k,i_{k+1}} &= {\mathbb P}_{\rm s}({\bs 1}(2)=i_2, {\bs 1}(k) = i_k, {\bs 1}(k+1) = i_{k+1}\,|\,{\bs 1}(1) = 0),\\
p^+_{i_2,i_k,i_{k+1}} &= {\mathbb P}_{\rm s}({\bs 1}(2)=i_2, {\bs 1}(k) = i_k, {\bs 1}(k+1) = i_{k+1}\,|\,{\bs 1}(1) = 1);
\end{align*}
here the subscript `s' denotes that the process is assumed to be in stationarity at time $1$ (which we have assumed throughout). The probabilities ${\bs p}^\pm:=(p^\pm_{0,0,0},p^\pm_{1,0,0},\ldots, p^\pm_{1,1,1})$ can be determined as in Remark \ref{R2}. 

Now condition on the value of $A_n(1)$; say it is $m\in\{0,1,\ldots,n\}$. Define 
\[A_n[i_2,i_{k},i_{k+1}] :=\# \big\{j\in\{1,\ldots,n\}:{\bs 1}_j(2)=i_2, {\bs 1}_j(k) = i_k, {\bs 1}_j(k+1) = i_{k+1}\big\}. \]
Then consider the random vector ${\bs A}_n\equiv (A_n[0,0,0], A_n[0,0,1],\ldots, A_n[1,1,1])\in\{0,\ldots,n\}^8$. It is evident that ${\bs A}_n$, conditional on $A_n(1)=m$, is the sum of two independent multinomially distributed quantities: in \michel{standard} notation,
\begin{equation}\label{eq:mult}{\bs A}_n\sim {\mathbb M}{\rm ult} \big(m, {\bs p}^+\big)+{\mathbb M}{\rm ult} \big(n-m, {\bs p}^-\big).\end{equation}
Now note that we can write, for a suitably chosen matrix $M$ of dimension $3\times 8$, 
\[\big(\begin{array}{c c c}
A_n(2)& A_n(k)& A_n(k+1)
\end{array}
\big)^\top = M {\bs A}_n; \]
for instance, $A_n(2) = A_n[1,0,0]+A_n[1,0,1]+ A_n[1,1,0]+ A_n[1,1,1]$. 
Then we write
\begin{align*}
   {\mathbb E}_{\rm s}\big[A_n(1)&A_n(2)A_n(k)A_n(k+1)\big]=\sum_{m=0}^n {\mathbb P}(A_n(1)=m)\,m\,
   {\mathbb E}_{\rm s}\big[A_n(2)A_n(k)A_n(k+1)\,|\,A_n(1)=m\big]\\
   &= \sum_{m=0}^n {\mathbb P}(A_n(1)=m)\,m\,{\mathbb E}_{\rm s}\big[\big( M {\bs A}_n\big)_1\big( M {\bs A}_n\big)_2\big( M {\bs A}_n\big)_3\,|\,A_n(1)=m\big].
\end{align*}
The last step is to plug in \eqref{eq:mult}, and to use that $A_n(1)\sim{\mathbb B}{\rm in}(n,\varrho).$
Then we end up with a weighted sum of mixed moments (of order up to 4) pertaining to multinomial distributions.

\vb

    \js{In the light of existing versions of the central limit theorem, it is anticipated that a necessary condition for asymptotic normality is the finiteness of \MM{$v_0$, $v_1$, and $c_{01}$.}
    Earlier in this section,} this finiteness was verified for the case of geometric on- and off-times, by explicitly identifying these three quantities. \michel{As demonstrated in Appendix \ref{AppB}, one can prove that, in case the on-times are ${\mathbb P}{\rm ar}(C_X,\alpha)$, and the off-times are ${\mathbb P}{\rm ar}(C_Y,\beta)$ for $\beta>\alpha$, one needs to have that $\min\{\alpha,\beta\}>2$, as stated in the following result.}

    \begin{proposition}\label{P3}
        \michel{Consider the case that $X\sim{\mathbb P}{\rm ar}(C_X,\alpha)$ and $Y\sim{\mathbb P}{\rm ar}(C_Y,\beta)$ with $C_X, C_Y>0$. If $\min\{\alpha,\beta\}>2$, then $v_0$, $v_1$, and $c_{01}$ are finite.}
    \end{proposition}

\begin{remark}
{\em
   Upon inspecting the derivation in Appendix \ref{AppB}, it can be seen that the conditions can be relaxed in various respects. In its most general form, one has to require that either ${\mathbb P}(X\geqslant i)$ or ${\mathbb P}(Y\geqslant i)$ has a regularly varying tail \cite{BGT} of index strictly larger than 2, and that the other tail is at most equally heavy. This for instance means that $v_0$, $v_1$, and $c_{01}$ are finite if $X\sim{\mathbb P}{\rm ar}(C,\alpha)$ with $\alpha>2$ and $Y\sim {\mathbb G}(q)$, or if $X\sim{\mathbb P}{\rm ar}(C,\alpha)$ with $\alpha>2$ and $Y\sim {\mathbb W}(\lambda,\beta)$.}\hfill$\Diamond$
   \end{remark}

\section{Beyond the number of edges: other subgraph counts}\label{SUB}
So far, we have set up a parametric inference procedure to estimate the on- and off-time distributions, based on observations of the evolution of the total number of edges. 
However, using essentially the same methodology, one can set up a similar procedure if the information available corresponds to other subgraph counts. In this section we demonstrate this for the case that we have access to observations of the number of triangles $T_N(k)$ at times $k\in\{1,\ldots,K\}$, where $N\in{\mathbb N}$ is the number of vertices of our dynamic Erd\H{o}s-R\'enyi graph, and we also present its counterpart for the case that we observe the number of wedges $W_N(k)$. 

\michel{Subgraph counts play a crucial role in random graph research. In this context we mention \cite{CHA}, a methodological breakthrough when it comes to the rare-event behavior of static Erd\H{o}s-R\'enyi random graphs. 
Various papers considered the likelihood of certain subgraph counts  {\it simultaneously} attaining specific values, sometimes giving rise to phase transitions. A few examples of papers in this spirit are \cite{STAR, KEN, RAD}. Regarding the relevance of subgraph counts in various application areas, see e.g.\ \cite{GAR}.}

First we notice that setting up likelihood-based methods would lead to serious complications, even in the case of geometrically distributed on- and off-times. To see this, observe that even for geometrically distributed on- and off-time, the process $\{T_N(k)\}_{k\in{\mathbb N}}$ is {\it not} Markovian, as one can have rather different configurations that lead to the same number of triangles. As we argue in the remainder of this section, one can develop a method of moments, though. 

Denote by ${\bs 1}_{i,j}(k)$ the indicator function that corresponds to the existence of the edge between the vertices $i$ and $j$ at observation $k\in\{1,\ldots,K\}$, for $i\not= j$ and $i,j\in\{1,\ldots,N\}$. The first moment equation is 
\[{\mathbb E}\,T_N(k) = {\mathbb E}\left[\sum_{i_1<i_2<i_3}{\bs 1}_{i_1,i_2}(k) \,
{\bs 1}_{i_1,i_3}(k)\, {\bs 1}_{i_2,i_3}(k) 
\right]=\binom{N}{3}\,\varrho^3.\]
Our next goal is to find a second moment equation. For $k\in\{1,\ldots,K-1\}$,
\begin{align*}
   {\mathbb E}\,T_N(k)&T_N(k+1) \\
   &= {\mathbb E}\left[\sum_{i_1<i_2<i_3}{\bs 1}_{i_1,i_2}(k) \,
{\bs 1}_{i_1,i_3}(k)\, {\bs 1}_{i_2,i_3}(k) \times \sum_{j_1<j_2<j_3}{\bs 1}_{j_1,j_2}(k+1) \,
{\bs 1}_{j_1,j_3}(k+1)\, {\bs 1}_{j_2,j_3}(k+1) 
\right].
\end{align*}
Let $a_m$, for $m=0,1,2,3$, the number of sets $\{j_1,j_2,j_3\}$, for a given set $\{i_1,i_2,i_3\}$, that have $m$ elements in common. It can be verified that
\[a_0=\binom{N-3}{3}=\binom{N}{3}-a_1-a_2-a_3,\:\:\:\:a_1=3\binom{N-3}{2},\:\:\:\: a_2 = 3(N-3),\:\:\:\:a_3 =1.\]

\begin{figure}
\centering
\subcaptionbox{\:{$m = 1$}}
{\includegraphics[width=0.35\linewidth]{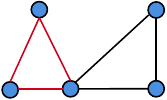}}
\qquad\qquad
\subcaptionbox{\:{$m = 2$}}
{\includegraphics[width=0.35\linewidth]{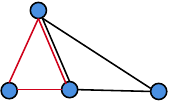}}
\caption{Triangle count, cases $m=1,2$. The {red and black} lines represent edges at time $k$ and $k+1$, respectively. The figure illustrates the elementary fact that if two triangles have one vertex in common ($m=1$), then they do not share an edge, and if they have two vertices in common ($m=2$), then they share one edge.} \label{fig:TriIllustarte}
\end{figure}
 Observe that (i)~if $\{i_1,i_2,i_3\}$ and $\{j_1,j_2,j_3\}$ have zero elements or one element in common, then they do not share an edge, (ii)~if they have two elements in common, then they share an edge, and~(iii)~if they are identical, then the corresponding triangles match (in that they share three edges); see Figure \ref{fig:TriIllustarte} for a pictorial illustration. For each instance corresponding to scenario (i), the probability of seeing a triangle at time $k$ as well as $k+1$ is therefore $\varrho^6$, for each instance corresponding to scenario (ii) it is $\varrho^3\cdot\varrho^2(1-\bar f_1)=\varrho^5(1-\bar f_1)$, and for each instance corresponding to scenario (iii) it is $\varrho^3(1-\bar f_1)^3.$
It thus follows that
\[{\mathbb E}\,T_N(k)T_N(k+1) = \binom{N}{3}\Big(
a_0 \varrho^6 + a_1 \varrho^6 + a_2 \varrho^5(1-\bar f_1)
+ a_3
\varrho^3(1-\bar f_1)^3\Big).\]
Combining arguments developed earlier in this paper, one can compute ${\mathbb E}\,T_N(k)T_N(k+\ell)$ for any $\ell\in{\mathbb N}$, so as to obtain additional moment equations (needed if the number of unknown parameters is larger than 2). 

A similar reasoning can be followed when setting up a method of moments for any type of subgraph count. One could, for instance, consider the number of `wedges' (i.e., groups of three connected vertices) 
\[W_N(k) = \sum_{i_1\not=i_2<i_3\not=i_1}{\bs 1}_{i_1,i_2}(k) \,
{\bs 1}_{i_1,i_3}(k);\]
here the index $i_1$ corresponds to the `center vertex', while $i_2$ and $i_3$ are the `end edges'.
The mean of the wedge count $W_N(k)$ is
\[{\mathbb E}\,W_N(k) = {\mathbb E}\left[\sum_{i_1\not=i_2<i_3\not=i_1}{\bs 1}_{i_1,i_2}(k) \,
{\bs 1}_{i_1,i_3}(k)
\right]=3 \, \binom{N}{3}\,\varrho^2,\]
whereas the corresponding second moment equation is
\[{\mathbb E}\,W_N(k)W_N(k+1) = \binom{N}{3}\Big(
9 a_0 \varrho^4 + 9 a_1 \varrho^4 + {5a_2}\varrho^4 + {4a_2} \varrho^3(1-\bar f_1)
+ 3a_3
\varrho^2(1-\bar f_1)^2 + 6a_3\varrho^3(1-\bar f_1)\Big).\]
\begin{figure}[!ht]
\centering
\subcaptionbox{\:{$m = 2$. The first line shows that if the center node at $k+1$ is not from the vertices that are present at time $k$, then the two wedges have no edges in common. The second line shows that if the center node at $k+1$ is one of the end nodes at $k$, then the two wedges have one edge in common if the center node at $k$ belongs to the vertices at $k+1$, and no edge in common otherwise. The third line shows that if the center node at $k+1$ is the center node at $k$, then there is one edge in common.}}
{\includegraphics[width=0.8\linewidth]{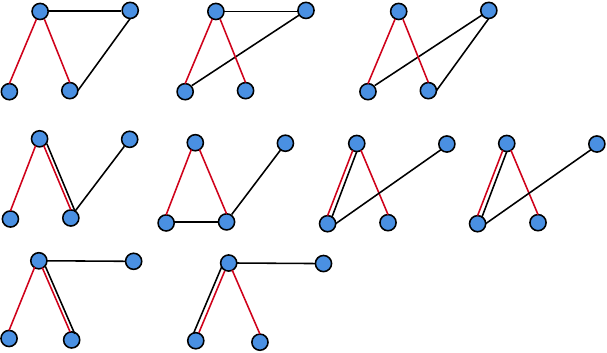}} 
\subcaptionbox{\:{$m = 3$. If the center node at $k+1$ is the same as the center node at $k$, then it is the same wedge; otherwise, there is only one edge in common.}}
{\includegraphics[width=0.8\linewidth]{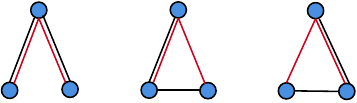}} 
\caption{{Wedge count, cases $m=2,3$. The {red and black} lines represent edges at time $k$ and $k+1$, respectively.}} \label{fig:WedgeIllustarte}
\end{figure}

This expression can be understood as follows. To form a wedge, three vertices are needed, and each vertex can be the center node at each of the times, so there are in total nine options. 
\begin{itemize}
\item[$\circ$]
Recall that $a_0$ and $a_1$ relate to the scenarios that $\{i_1,i_2,i_3\}$ and $\{j_1,j_2,j_3\}$ have zero elements, respectively one element in common. 
    For the expressions related to $a_0$ or $a_1$ this means that in these scenarios do not correspond to an edge being present at time $k$ as well as $k+1$. 
    \item[$\circ$]
    To understand the terms related to $a_2$, we distinguish two cases. One case relates to the configurations in which the two common vertices are both non-center at time $k$, whereas in the other case one of them is the center vertex at time $k$. The first case can occur in three ways, where it is directly verified (see the first line of Figure~\ref{fig:WedgeIllustarte}({\sc a})) that one cannot have that there is an edge present at time $k$ as well as $k+1$. 
    The second case can occur in six ways, with five of them resulting in no edges in common and one in one edge in common (see the second and third lines of Figure~\ref{fig:WedgeIllustarte}({\sc a})). Putting the above together, in four out of nine cases there are no edges in common (leading  to a contribution $\varrho^4$) and in five out of nine cases there is one edge in common (leading to a contribution $\varrho^3(1-\bar f_1)$). 
    \item[$\circ$]A similar reasoning applies to the terms related to $a_3$: if the center vertex at $k$ is the same as at $k+1$, the edges are the same at both times. As depicted in the first subfigure of Figure~\ref{fig:WedgeIllustarte}({\sc b}), this gives three options, each of them leading to a contribution $\varrho^2(1-\bar f_1)^2$). Otherwise, only one edge is in common. As depicted in the last two subfigures of Figure~\ref{fig:WedgeIllustarte}({\sc b}), this gives six options, each of them leading to a contribution $\varrho^3(1-\bar f_1)$).
\end{itemize}

\begin{remark}\label{R3}{\em 
In principle any other subgraph count can be handled analogously, i.e., by properly counting. One could also think of estimators based on {\it multiple} subgraph counts (e.g.\ squares and wedges). As these extensions in principle follow the same line of reasoning, we have decided to leave them out here.} $\hfill\Diamond$
\end{remark}

\section{Numerical Experiments}\label{NUM}
In this section we present numerical experiments, by which we empirically assess the performance of our estimation procedure. 

\medskip

In the first series of experiments, we estimate the parameters of the on- and off-times distributions based on the number of edges. We consider the following four instances:
\begin{itemize}
    \item[$\circ$] $X\sim {\mathbb G}(p)$ with $p=0.3$, and $Y\sim {\mathbb G}(q)$ with $q=0.8$.
    \item[$\circ$] $X\sim {\mathbb P}{\rm ar}(1,\alpha)$ with $\alpha = 3$, and $Y\sim {\mathbb P}{\rm ar}(1,\beta)$ with $\beta = 2.5$. 
    \item[$\circ$] $X\sim {\mathbb W}(1,\alpha)$ with $\alpha=0.5$, and $Y\sim {\mathbb G}(q)$ with $q=0.7$.
    \item[$\circ$] $X\sim {\mathbb P}{\rm ar}(C,\alpha)$ with $C=2$ and $\alpha=4$, and $Y\sim {\mathbb G}(q)$  with $q=0.7$.
\end{itemize}
In the first three instances there are two parameters to be estimated, and in the fourth instance three parameters. 
For each of these four settings the estimators were provided in Subsection~\ref{SS:EST}.

\begin{figure}
\centering
{${\mathbb G}(0.3)/{\mathbb G}(0.8)$} \\
\subcaptionbox{\footnotesize {$\overline{\hat p}^{(L)}=0.3000$ \, ($\sigma_p^{(L)}=0.0009)$}}
{\includegraphics[width=0.38401\linewidth]{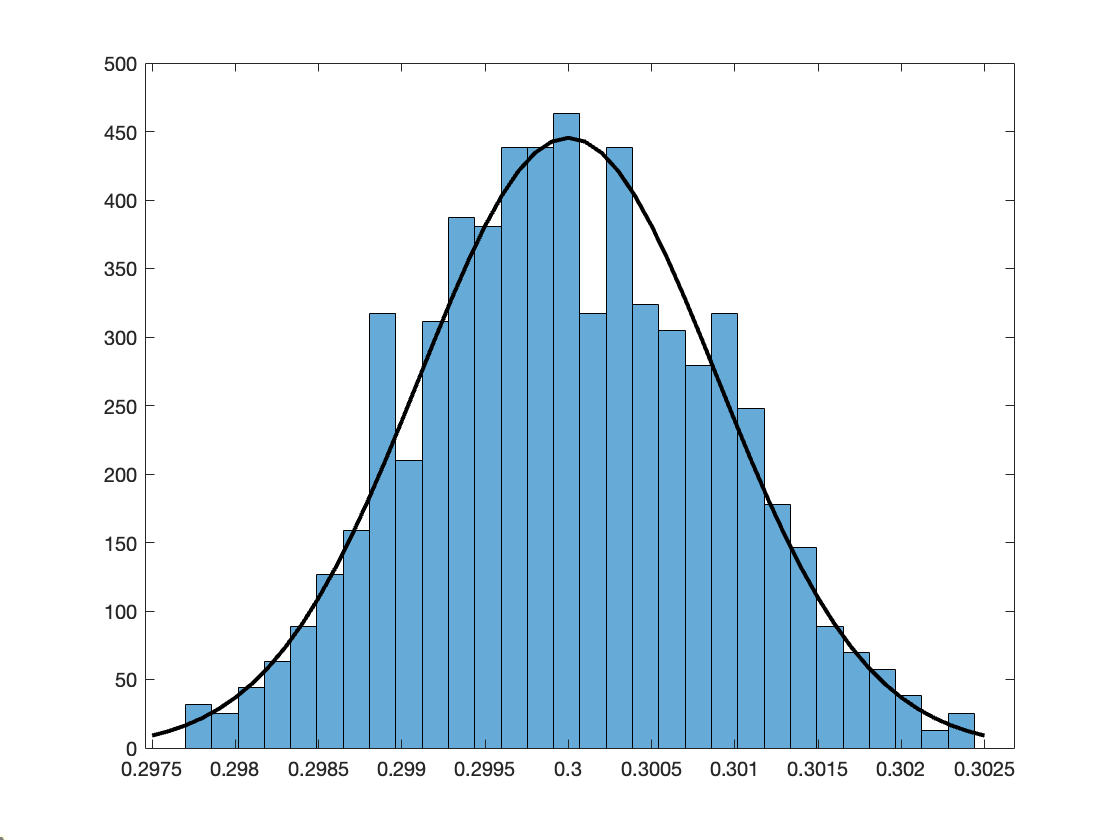}}
\subcaptionbox{\footnotesize {$\overline{\hat q}^{(L)} =0.8000$ \, ($\sigma_q^{(L)}=0.0024$)}}
{\includegraphics[width=0.38401\linewidth]{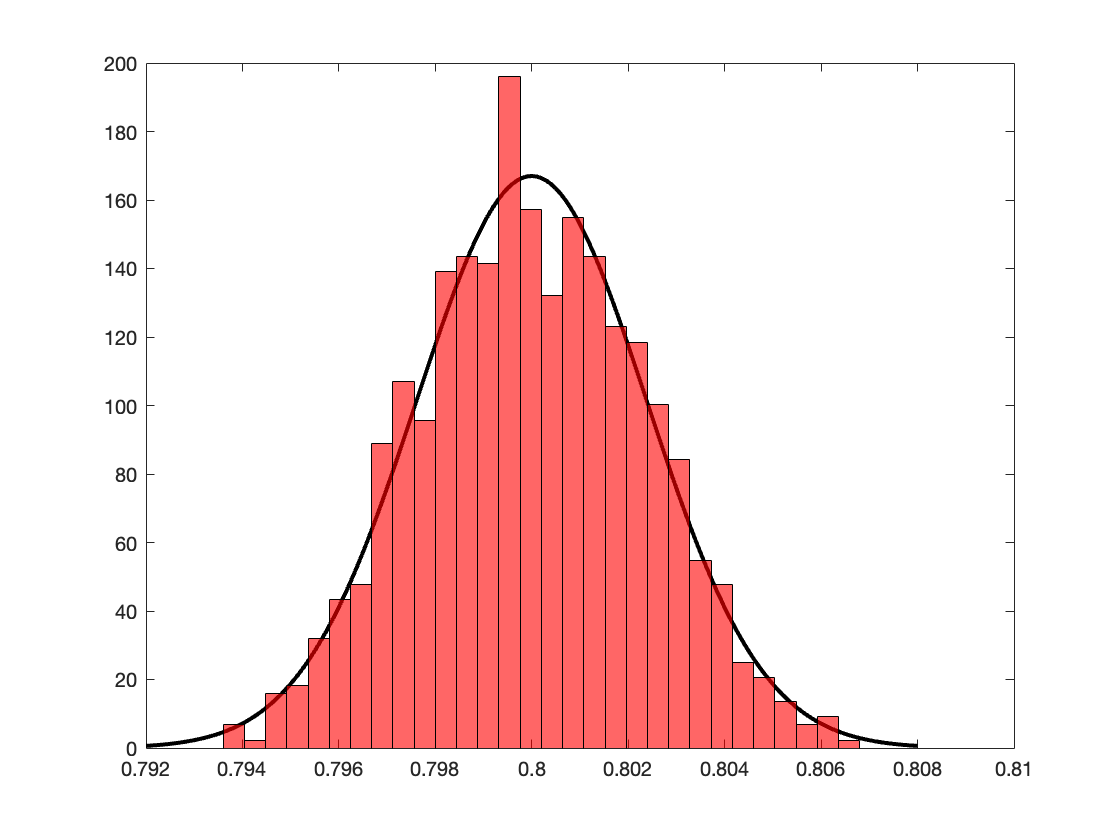}}\\
\vspace{2mm}
{${\mathbb P}{\rm ar}(1,3)/{\mathbb P}{\rm ar}(1,2.5)$}\\
\subcaptionbox{\michel{\footnotesize $\overline{\hat\alpha}^{(L)}=3.0002$ \, ($\sigma_\alpha^{(L)}=0.0269$)}}{\includegraphics[width=0.38401\linewidth]{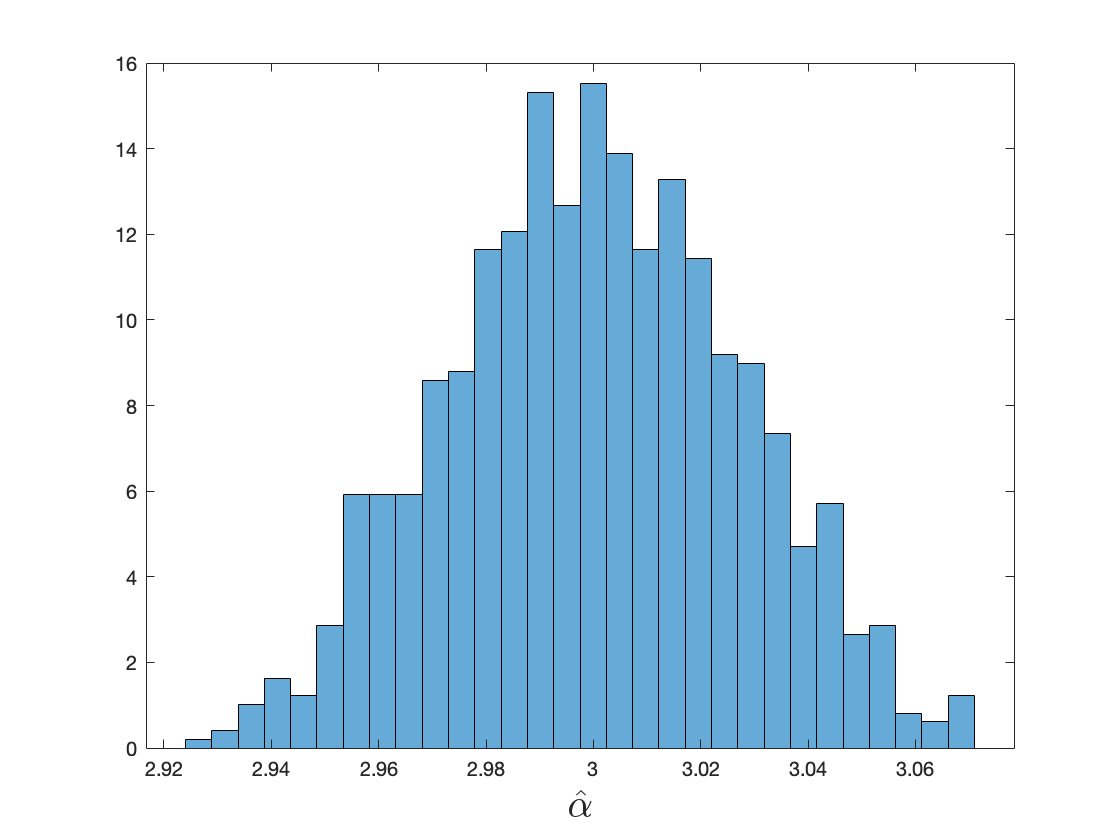}}
\subcaptionbox{\michel{\footnotesize $\overline{\hat\beta}^{(L)}=2.4999$ \, ($\sigma_\beta^{(L)}=0.0153$)}}
{\includegraphics[width=0.38401\linewidth]{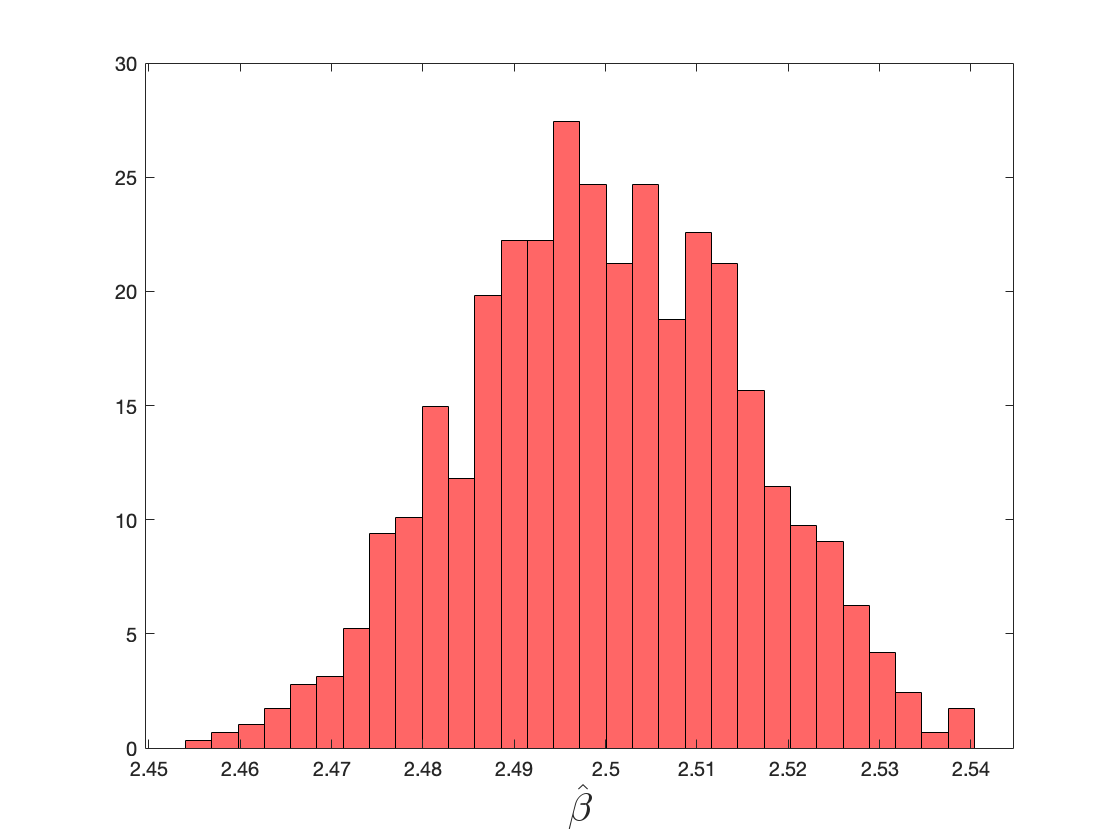}}\\
\vspace{2mm}
\michel{${\mathbb W}(1,0.5)/{\mathbb G}(0.7)$}\\
\subcaptionbox{\michel{\footnotesize  $\overline{\hat\alpha}^{(L)}= 0.5000$ \, ($\sigma_\alpha^{(L)}=0.0014$)}}
{\includegraphics[width=0.38401\linewidth]{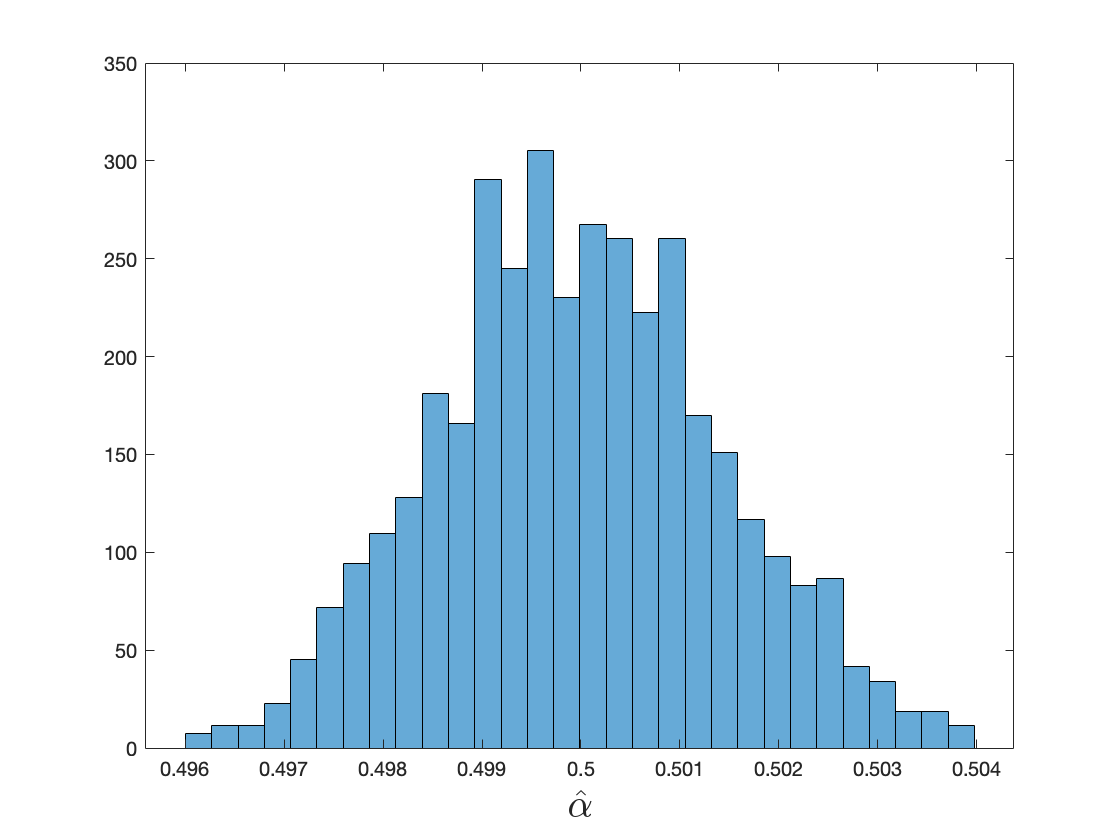}}
\subcaptionbox{\michel{\footnotesize $\overline{\hat q}^{(L)} =0.7000$ \, ($\sigma_q^{(L)}=0.0028$)}}
{\includegraphics[width=0.38401\linewidth]{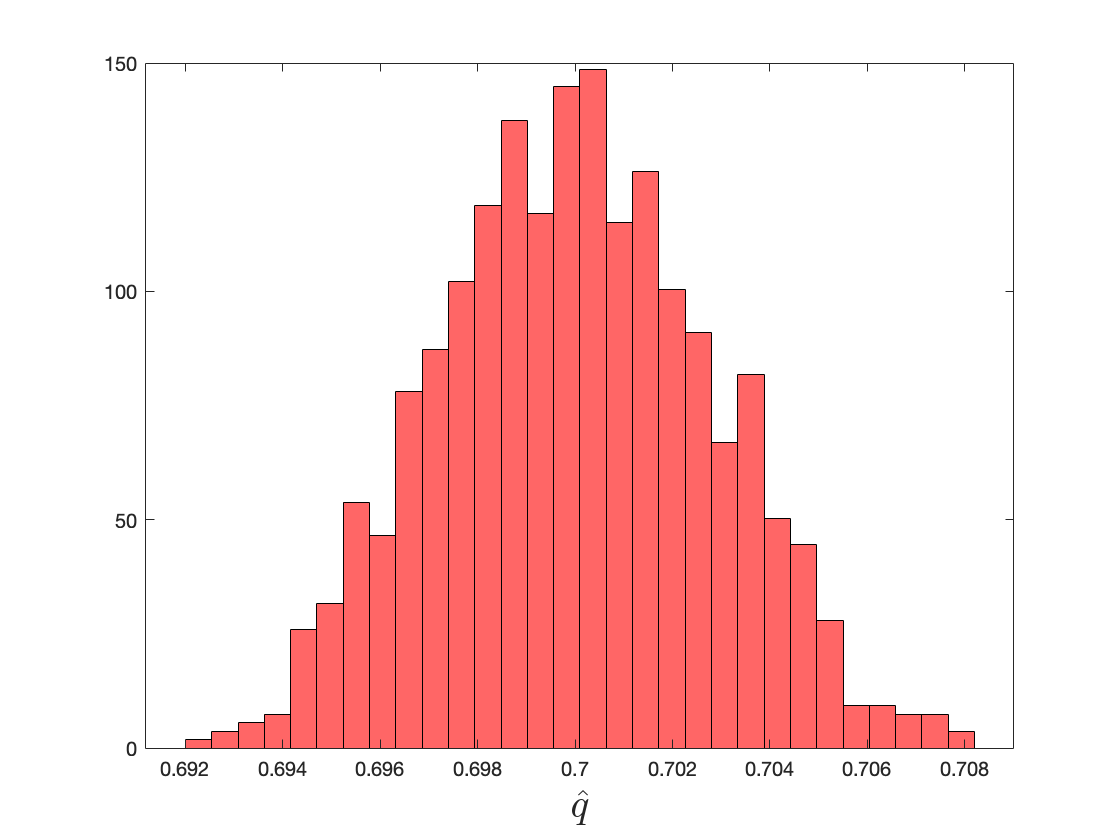}}\\
\vspace{2mm}
\michel{${\mathbb P}{\rm ar}(2,4)/{\mathbb G}(0.7)$}\\
\subcaptionbox{\michel{\footnotesize $\overline{\hat C}^{(L)}=2.0259$ \,($\sigma_C^{(L)}=0.3239$),\\\hspace*{6mm}$\overline{\hat \alpha}^{(L)}= 4.0386$ \,($\sigma_\alpha^{(L)}=0.4705$)}}
{\includegraphics[width=0.38401\linewidth]{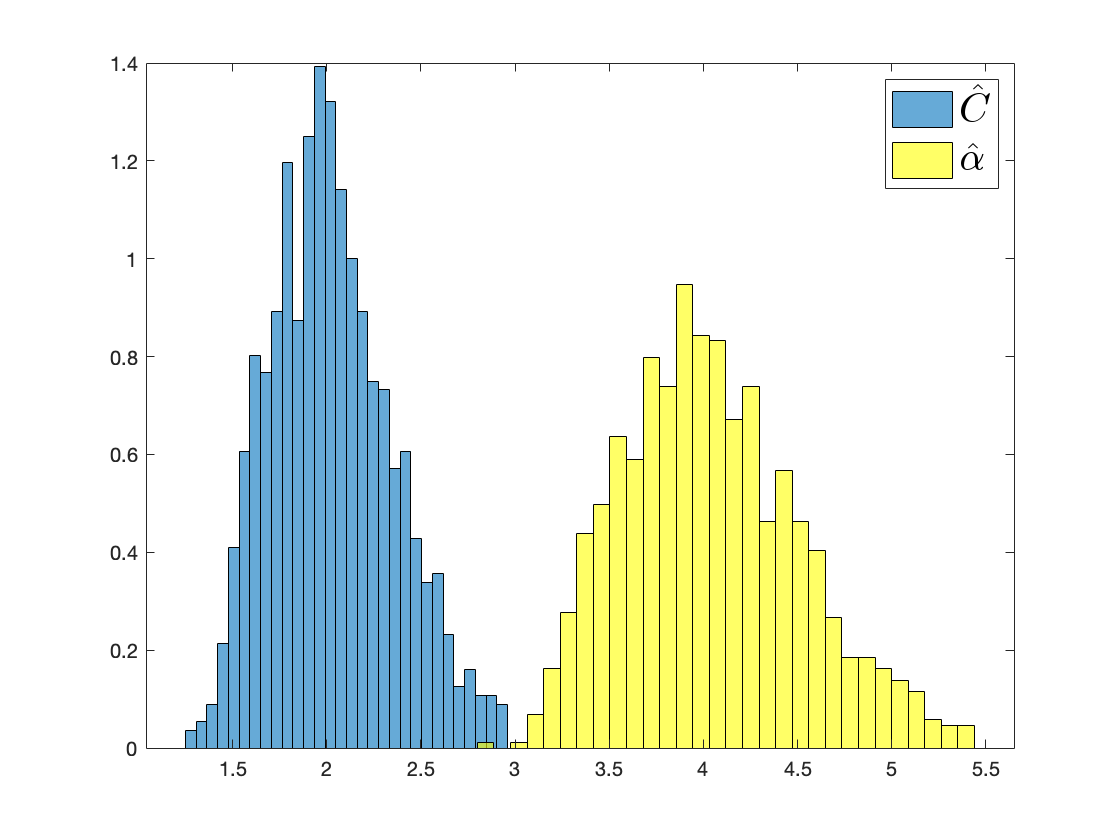}}
\subcaptionbox{\michel{\footnotesize $\overline{\hat q}^{(L)} =0.7000$ \, ($\sigma_q^{(L)}=0.0078$)}}
{\includegraphics[width=0.38401\linewidth]{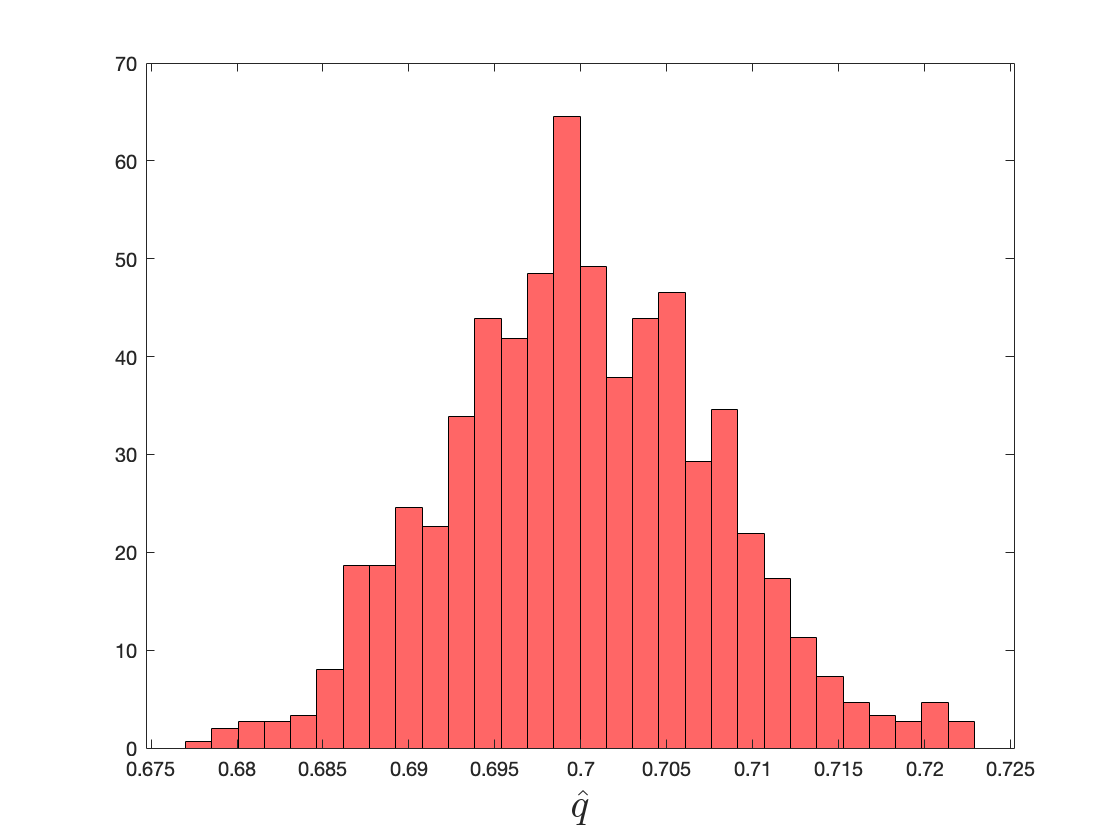}}
\caption{Estimation from the method of moments when $n=100$ and $K = 10^5$ (results from $L=1000$ experiments). Standard deviations are provided in parentheses.} \label{fig:momhist}
\end{figure}

\michel{We have generated traces of $\{A_n(k)\}$ for $k=1,\ldots,10^5$ and $n=100$, and estimated all parameters. \MM{We performed this procedure $L=1000$ times for each instance; each of the four rows in  Figure~\ref{fig:momhist} corresponds to one of the four instances defined above.
For example, the two panels in first row of Figure~\ref{fig:momhist}, which represent the case where $X\sim {\mathbb G}(p)$ and $Y\sim {\mathbb G}(q)$, display histograms pertaining to the estimates \[\hat p^{(1)}_{n,K},\ldots, \hat p^{(L)}_{n,K},\quad \hat q^{(1)}_{n,K},\ldots, \hat q^{(L)}_{n,K}.\] 
In the sequel we denote by \begin{align}\label{pqdef}
\overline{\hat p}^{(L)}_{n,K}& \equiv \overline{\hat p}^{(L)}:= \frac{1}{L}\sum_{\ell=1}^L \hat p^{(\ell)}_{n,K},\\ \label{pqdef2}
\overline{\hat q}^{(L)}_{n,K}& \equiv \overline{\hat q}^{(L)}:= \frac{1}{L}\sum_{\ell=1}^L \hat q^{(\ell)}_{n,K}\end{align} the means of these vectors, and by $\hat\sigma_{p,n,K}^{(L)}\equiv \hat\sigma_p^{(L)}$ and $\hat\sigma_{q,n,K}^{(L)}\equiv \hat\sigma_q^{(L)}$ their standard deviations. The histograms shown in the remaining three rows of Figure \ref{fig:momhist} should be interpreted in the same manner.}}

\medskip

\js{In the experiments reported, the parameters are typically estimated accurately. 
The exception is the setting of Figure~\ref{fig:momhist}.(H) where the confidence intervals are wider, potentially a consequence of the fact that in that experiment {\it three} (rather than two) parameters are estimated. This means that in that case the moment equations contain $s_{n,0}$, $s_{n,1}$, and $s_{n,2}$ (rather than just $s_{n,0}$ and $s_{n,1}$), which could have led to the loss of accuracy discussed in Remark \ref{RR1}.  In additional experiments involving Pareto-type distributions (with a fixed number of unknown parameters), not reported in this paper, we have observed that standard deviations increase when the distribution's shape parameter approaches $2$ from above; this behavior is confirmed by the reasoning in Appendix \ref{AppB}.}

\medskip

\begin{figure}
    \centering
    {\includegraphics[width=0.48\linewidth]{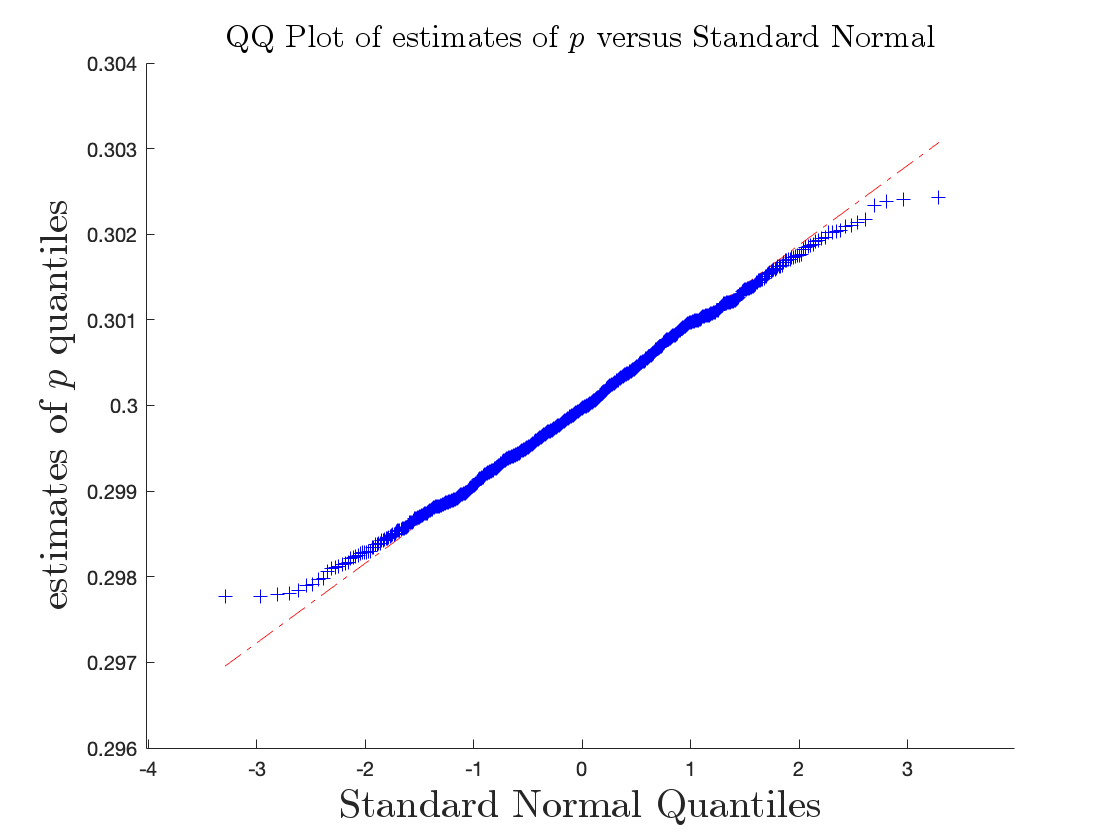}}
    {\includegraphics[width=0.48\linewidth]{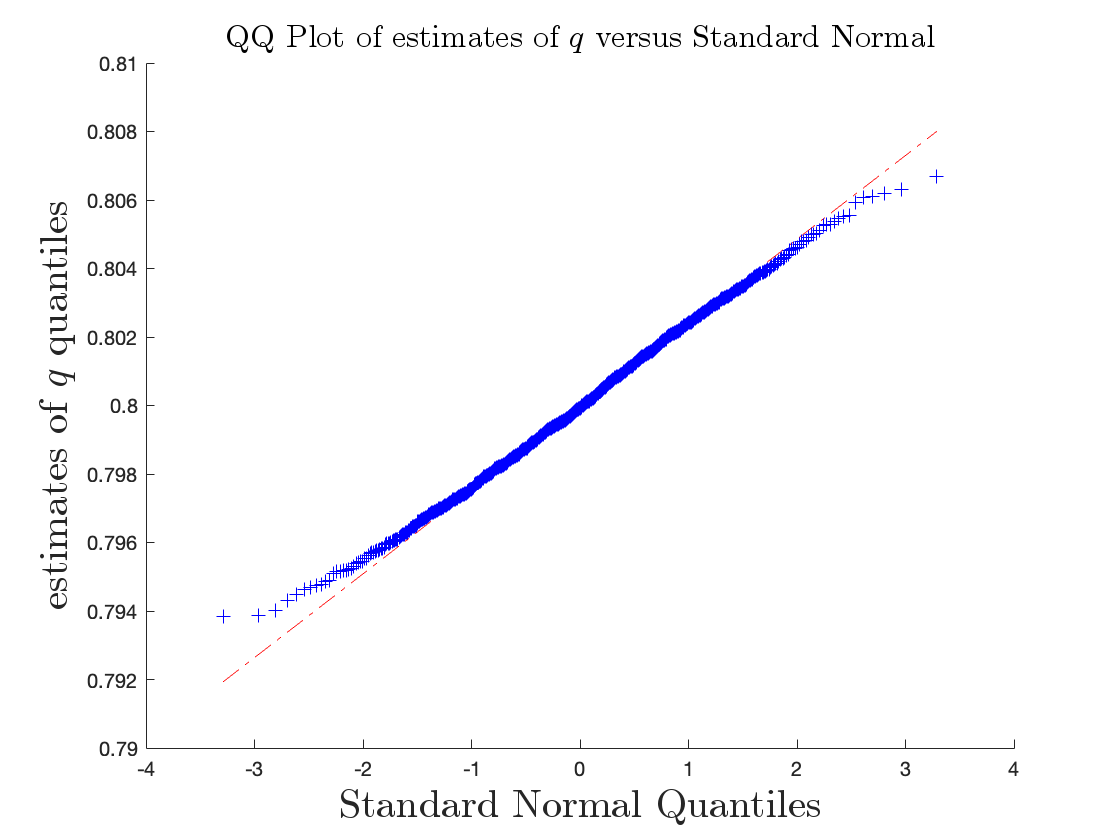}}
    \caption{\label{FIGQQ}Q-Q plots of estimates $\overline{\hat p}^{(L)}$ and $\overline{\hat q}^{(L)}$ in the ${\mathbb G}(0.3)/{\mathbb G}(0.8)$ case. We picked $n = 100$  and $K =10^5$ (results from $L=1000$ experiments).}
    \label{fig:QQPlot}
\end{figure}

We in addition generated QQ-plots to assess the asymptotic normality. All of these plots showed nearly straight lines, as desired; see for example Figure \ref{FIGQQ} for the ones corresponding to the case of geometric on- and off times. Only in the `far tails' we observe a significant deviation from the straight line (which will disappear when choosing $K$ larger).

\medskip

\begin{figure}
\centering
{\vspace*{1cm}Estimations based on triangle counts} \\
\subcaptionbox{\footnotesize $\overline{\hat p}^{(L)}=0.3001 \,\,\, (\sigma^{(L)}_p=0.0029)$}
{\includegraphics[width=0.401\linewidth]{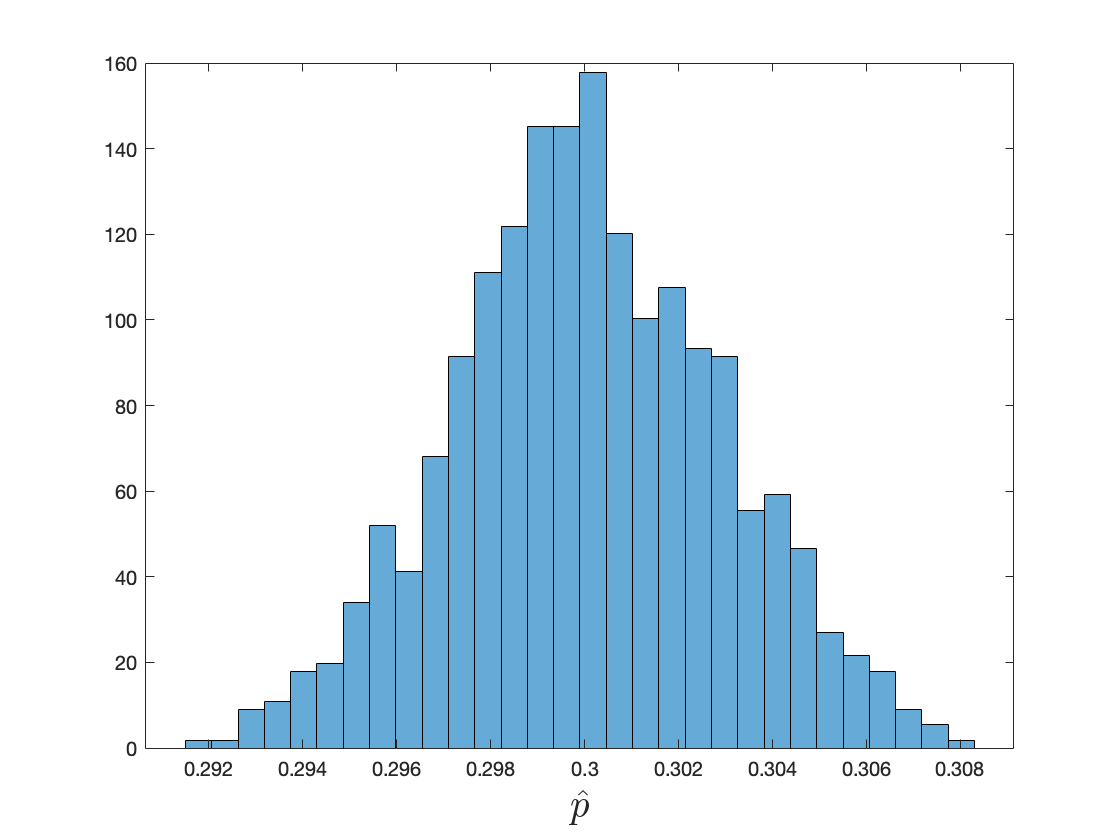}}
\subcaptionbox{\footnotesize $\overline{\hat q}^{(L)} =0.8002 \,\,\, (\sigma_q^{(L)}=0.0076)$}
{\includegraphics[width=0.401\linewidth]{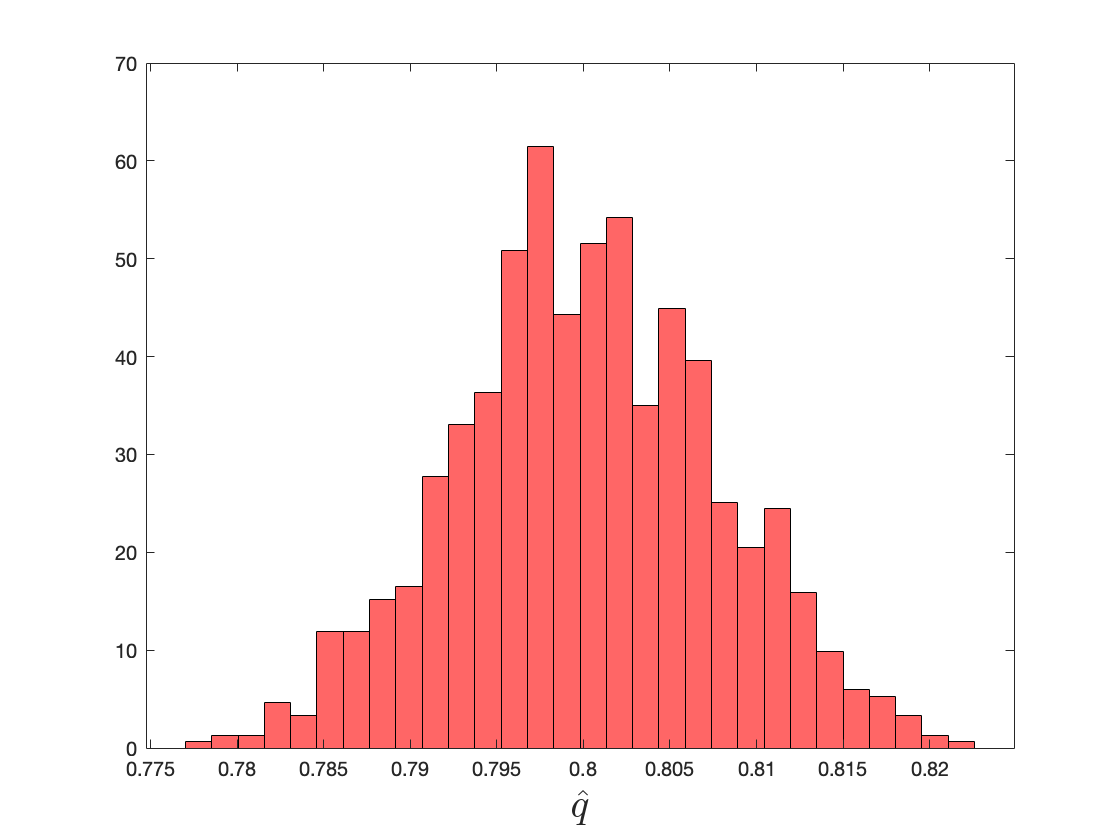}}\\
\vspace{4mm}
{Estimations based on wedge counts}\\
\subcaptionbox{\footnotesize $\overline{\hat p}^{(L)}= 0.3000 \,\,\,  (\sigma^{(L)}_p=0.0027)$}
{\includegraphics[width=0.401\linewidth]{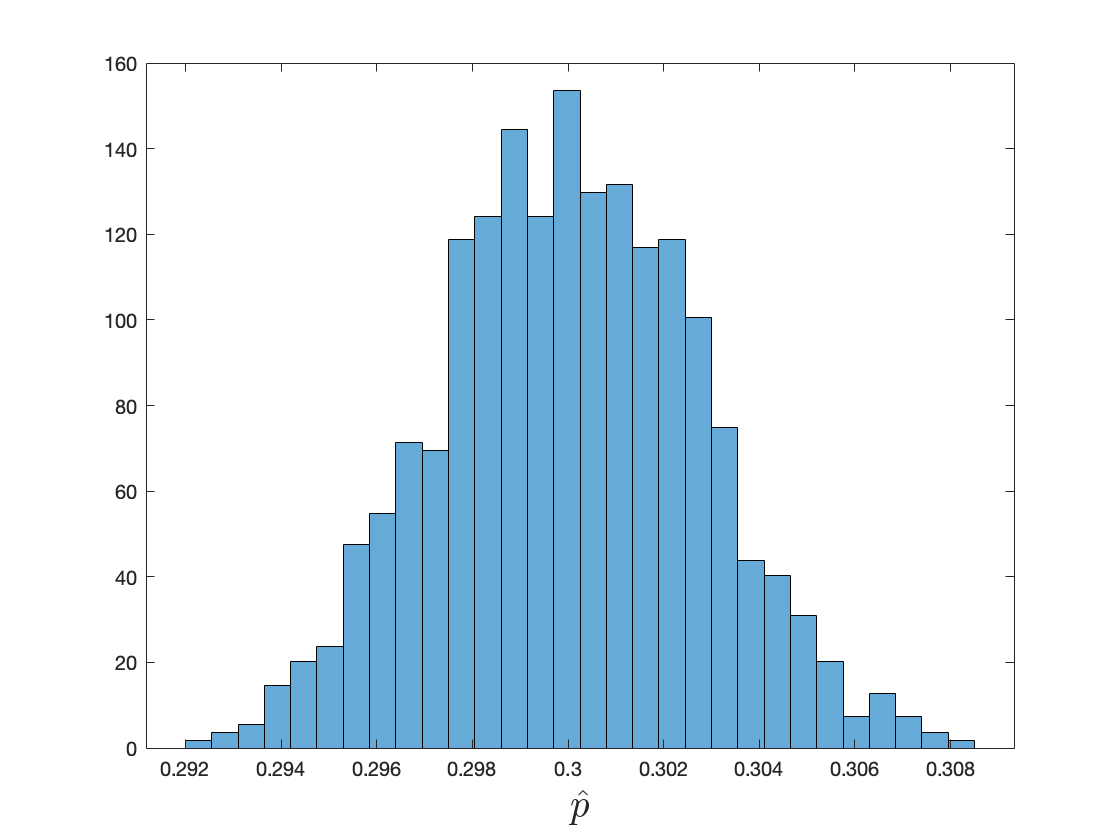}}
\subcaptionbox{\footnotesize $\overline{\hat q}^{(L)}=  0.8001 \,\,\,  (\sigma^{(L)}_q=0.0072)$}
{\includegraphics[width=0.401\linewidth]{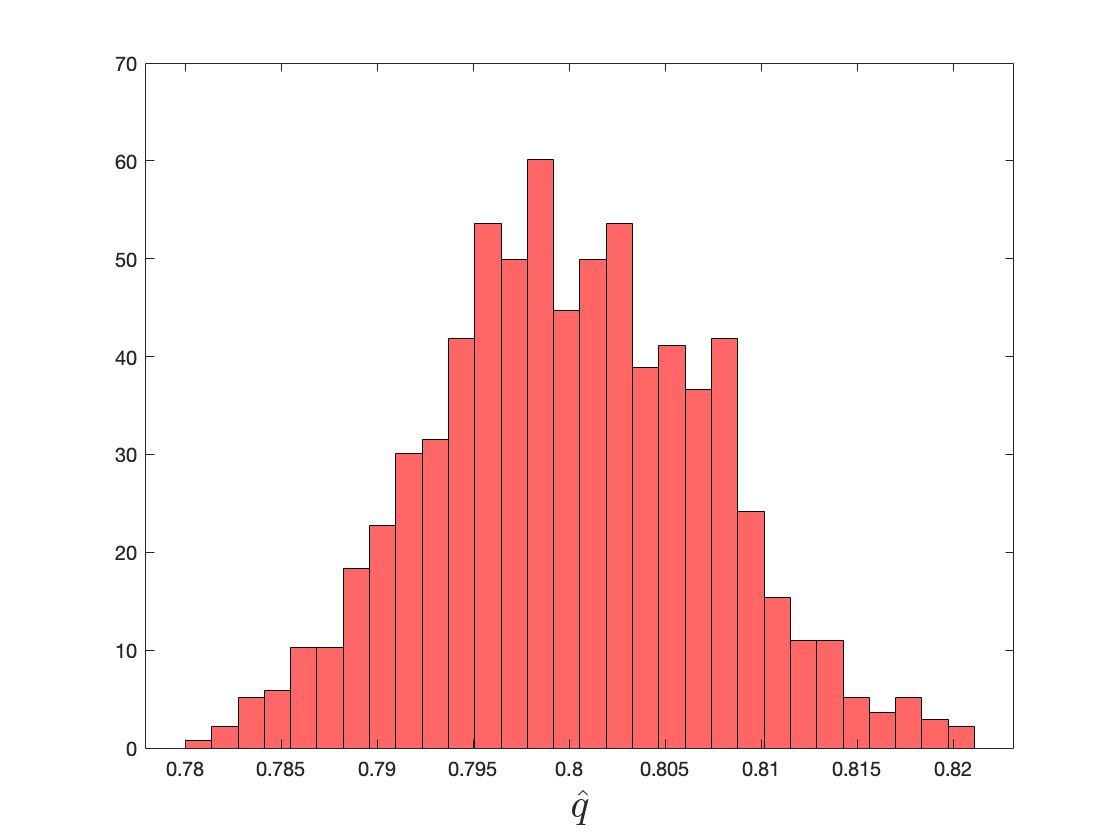}}
\caption{Estimation from the method of moments based on subgraph counts when $N = 20$ (i.e., $n = 190$) and $K = 10^4$ (results from $L=1000$ experiments). Standard deviations are provided in parentheses.} \label{fig:momTriWedgeHist}
\end{figure}

\jiesen{\MM{For the ${\mathbb G}(0.3)/{\mathbb G}(0.8)$ case in Figure \ref{fig:momhist}, we also plotted the probability density functions of a Normal distribution in which 
\begin{itemize}
    \item[$\circ$] the respective means are $\overline{\hat p}^{(L)}$ and $\overline{\hat q}^{(L)}$, as given in \eqref{pqdef}--\eqref{pqdef2},
    \item[$\circ$] the respective variances as given in Proposition \ref{P1}, with $p$ and $q$ replaced by  $\overline{\hat p}^{(L)}$ and $\overline{\hat q}^{(L)}$.
\end{itemize}
These density functions are in good agreement with the corresponding histograms.} 

The code to generate the plots is in the {\small \tt checkVariancePQ.m} file in the folder {\small \tt Asymptotic} in 
{\small \url{https://github.com/encwang/OnOffEstimation}}.
All the plots in Figures \ref{fig:momhist} and \ref{fig:momTriWedgeHist} can be obtained by running the {\small \tt Plot.m} file, and the data used in generating the plots in this paper is stored in the {\small \tt .mat} file, in each of the folders.
}


\medskip

{In another series of experiments we performed estimations based on triangle and wedge counts. We focused on the case ${\mathbb G}(0.3)/{\mathbb G}(0.8)$; we set the number of vertices $N = 20$, so that $n = N(N-1)/2 = 190$.
In each run we record $\{A_n(k)\}$ for $k = 1,\ldots, 10^4$, and count the number of triangles (wedges) at each $k$. The estimation procedure was performed \MM{$L=10^3$} times; estimates from these $10^3$ runs have been visualized via histograms in Figure \ref{fig:momTriWedgeHist}. The graphs confirm that the parameters are estimated accurately.}

\section{Discussion and concluding remarks}\label{DISC}

    In this paper we have proposed a procedure to estimate on- and off-time distributions in a dynamic Erd\H{o}s-R\'enyi graph,  only observing the aggregate number of edges. We developed an estimator based on the method of moments, which under mild conditions is asymptotically normal. 
    We also pointed out how the estimation procedure can be adapted if alternative subgraph counts are observed.

    \vb

    \michel{
    Our experiments show that our estimation procedure performs well, also in the context of heavy-tailed distributions (such as our Pareto-type distribution). This may seem remarkable, given that in our estimators we only use empirical information pertaining to the {\it lower tail} of the distributions of $X$ and $Y$. 
    In this context it should be realized that the favorable performance of our estimators is a consequence of the fact that we work in a {\it parametric} setting: we are given the parametric form of the probability mass functions (or, equivalently, the cumulative distribution functions) of the on- and off-times. 
    
    The analysis becomes intrinsically more involved when we are only given the {\it asymptotic} behavior of the on- and off-time distributions. A particularly relevant case is that in which  $X$ and/or $Y$ are  {asymptotically} of Pareto type. For instance, an objective could be to estimate $C$ and $\alpha$, given that we know that
    \[{\mathbb P}(X\geqslant i) \,i^\alpha \to C,\quad i\to\infty.\]
   \js{In such settings, one may consider setting up an estimator that incorporates information on correlations over longer timescales.}
    }

    \vb

We conclude this paper by discussing four possible directions for future research: (i) a first in which the on- and off-times are continuous non-negative random variables and in which the number of edges is observed \michel{at the times of a Poisson
process}, (ii) a second in which no features of the dynamic graph is observed, but rather features of a process evolving on it, (iii) a third in which the graph dynamics are affected by an unobservable background process, and (iv) a fourth relating to the regime that the random graph under consideration grows large.  

(i)~Suppose $X$ has a continuous distribution on ${\mathbb R}_+$ with density $f(\cdot)$ and Laplace-Stieltjes transform ${\mathscr F}(\cdot)$, and $Y$ a continuous distribution on ${\mathbb R}_+$ with density $g(\cdot)$ and Laplace-Stieltjes transform ${\mathscr G}(\cdot)$. Suppose that the number of edges is observed \michel{at the times of a Poisson
process}, \js{so that the inter-observation times are exponentially distributed, say with parameter $\xi>0$}. When starting the process in stationarity, \js{every edge is on with probability $\varrho$ and off with probability $1-\varrho$. The remaining on- and off- times have densities}
\[\bar f(x):= \frac{\int_{y=x}^\infty f(y)\,{\rm d}y}{{\mathbb E}\,X},\:\:\:\:\bar g(x):= \frac{\int_{y=x}^\infty g(y)\,{\rm d}y}{{\mathbb E}\,Y},\]
respectively. 

This model can be estimated in the same way as the discrete-time model featuring in this paper. We sketch the main steps. 
The observation $A_n(k)$ now corresponds to the number of edges at the $k$-th Poisson inspection epoch. 
Let $p_\xi(++)$ be the probability that the edge is on at an exponentially distributed time (with mean $\xi^{-1}$) given a fresh on-time started at time $0$, and let $p_\xi(+-)$, $p_\xi(-+)$, and $p_\xi(--)$ be defined analogously. Then, relying on the memoryless property of the exponential distribution,
\begin{align*}
  p_\xi(++) &= \int_0^\infty \int_0^t f(x) \,\xi e^{-\xi t} \,  p_\xi(-+) \,{\rm d}x\,{\rm d}t + \int_0^\infty \int_t^\infty f(x) \,\xi e^{-\xi t} \, {\rm d}x\,{\rm d}t = 1 - p_\xi(--)\,{\mathscr F}(\xi),
\end{align*}
where in the second equality we have swapped the order of integration and we have used $p_\xi(-+) = 1-p_\xi(--)$.
Likewise, 
\[ p_\xi(--) =  1 - p_\xi(++)\,{\mathscr G}(\xi),\]
so that
\[p_\xi(++) =\frac{1-{\mathscr F}(\xi)}{1-{\mathscr F}(\xi)\,{\mathscr G}(\xi)},\:\:\:\:p_\xi(--) =\frac{1-{\mathscr G}(\xi)}{1-{\mathscr F}(\xi)\,{\mathscr G}(\xi)}. \]
With $p^{\rm (res)}_\xi(++)$ the counterpart of $p_\xi(++)$ but then starting with a {\it residual} on-time (rather than a fresh one), we \js{similarly} obtain
\[p^{\rm (res)}_\xi(++) = 1 - p_\xi(--)\,{\mathscr F}^{\rm (res)}(\xi),\:\:\:\:
{\mathscr F}^{\rm (res)}(\xi) =\int_0^\infty e^{-\xi x}\bar f(x)\,{\rm d}x=\frac{1-{\mathscr F}(\xi)}{\xi\,{\mathbb E}X}.\]
Supposing the on- and off-times are each characterized by a single parameter, two moment equations are needed. These could be
\[{\mathbb E}A_n(k) = n\varrho,\:\:\:{\mathbb E}\big[A_n(k)A_n(k+1)\big]=n\varrho\,p^{\rm (res)}_\xi(++)+(n^2-n) \varrho^2.\]
In case there are more than two parameters, additional moment conditions need to be identified. This procedure is fully analogous to the one relied upon in the discrete-time case: for instance when there are three parameters, one needs to in addition evaluate (in self-evident notation) the probabilities $p^{\rm (res)}_\xi(+\hspace{-0.15em}+\hspace{-0.15em}+)$ and $p^{\rm (res)}_\xi(+\hspace{-0.15em}-\hspace{-0.15em}+)$, which can be done as above.

The setting in which the number of edges is observed at {\it equidistant} points in time, rather than \michel{at the times of a Poisson
process}, is considerably more challenging.

\medskip

(ii)~In a second class of related problems one does not observe features of the dynamic graph itself, but rather a process evolving on it. An example could be the model in which there is a fixed number of `random walkers', who (at every unit of time) either stay at the node where they reside or jump to one of their neighbors; a potential mechanism could be that if a vertex has $\ell$ neighbors, with probability $\ell /(\ell+1)$ one of these neighbors is selected (say in a uniform way) and with probability $1/(\ell+1)$ the random walker stays at the same vertex. 

\medskip

(iii)~A third class of problems relates to the situation in which the graph dynamics are affected by a background process. The most elementary example is that in which there is a Markovian background process on $\{1,2\}$, such that if this process is in state $i$ the probability of an existing edge turning off is $p_i$ and that of a non-existing edge turning on is $q_i$ (the two-dimensional modulated counterpart of the model with ${\mathbb G}(p)$ on-times and  ${\mathbb G}(q)$ off-times, that is). A complication is that the background process is not observed. In other contexts, similar problems have been tackled in e.g.\ \cite{GUN,GUN2, OKA}. For instance, in \cite{GUN2, OKA} increments of a Markovian arrival process are observed, from which both the parameters of the background process and the (state-dependent) arrival rates are estimated.

\medskip

    \michel{(iv)~Many results in the random graph literature relate to regimes in which the underlying graph becomes large. This can be done by sending $n$ to $\infty$, potentially in combination with some scaling of the random variables $X$ and $Y$. In our work, however, $n$ is held fixed; we consider properties of estimators when the number of observations $K$ of subgraph counts grows large. A potential topic for future research could concern an asymptotic regime in which $n$ and $K$ {\it jointly} grow. }

    \vb

\appendix

\section{Variances and covariance for geometric on- and off-times}\label{AppA}
\michel{In this appendix we compute $v_1$ and $c_{01}$ for the case of geometric on- and off-times. 
We start by computing $v_1$, which we do in four steps.

\medskip

{\it Step 1: decompose \js{$v_1$}.} 
Using the same reasoning as in the $v_0$ case, we need to evaluate
\begin{align*}\MRHM{v_1}=t_1+2 \sum_{k=2}^{\infty}t_k, \end{align*}
with $t_k:={\mathbb C}{\rm ov}\left(A_n(1)A_n(2),A_n(k)A_n(k+1)\right)$.
The cases $k=1$, and $k\in\{2,3\ldots\}$ have to be treated separately.

\medskip

{\it Step 2: calculation of stationary moments.} To facilitate the computations, we introduce a number of auxiliary objects, relating to the marginal distribution of $A_n(1)$ and the joint distribution of $A_n(1)$ and $A_n(2)$. Let 
\[m_i:={\mathbb E}\big[\big(A_n(1)\big)^i\big] \qquad \qquad \mathfrak{m}_i:={\mathbb E} \big[A_n(1)\big(A_n(2)\big)^i\big] .\]
Because $A_n(1)\sim {\mathbb B}{\rm in}(n,\varrho)$, the moments $m_i$ can be derived in the standard manner, e.g., by differentiation of the moment generating function. \js{It is readily checked that}
\begin{align*}
    & m_1 = n\varrho ,\\
    & m_2 = n(n-1)\varrho^2 + n\varrho ,\\
    & m_3 =  n(n-1)(n-2) \varrho^3 + 3n(n-1)\varrho^2 + n \varrho, \\
    & m_4 = n(n-1)(n-2)(n-3) \varrho^4 + 6n(n-1)(n-2)\varrho^3 + 7n(n-1)\varrho^2 + n\varrho \,.
\end{align*}
Due to 
 \[A_n(2) \sim {\mathbb B}{\rm in}(A_n(1), 1-p) + {\mathbb B}{\rm in}(n-A_n(1),q),\]
 with the usual independence between the two quantities on the right-hand side, 
we find after elementary calculations that
\js{
\begin{align}
    & \mathfrak{m}_1 = \lambda m_2+n q \,m_1, \label{eq:m1}\\
    & \mathfrak{m}_2 = \lambda^2 m_3 + \lambda (p-q+2nq) m_2 + (n q (1-q)+n^2 q^2) m_1, \label{eq:m2}\\
    & \mathfrak{m}_3 = \lambda^3 m_4 + 3 \lambda^2 (p-q + n q) m_3 \label{eq:m3}\\
    & \qquad \qquad\:\:\:\:\:+ \MM{\lambda}\, (3 n^2 q^2+3 n p q-6 n q^2 + 3 n q + 2 p^2 - 2 p q - p + 2 q^2 - q) m_2 \notag \\
    & \qquad \qquad\:\:\:\:\:+ (n(n - 1)(n - 2) q^3 + 3 n (n - 1) q^2 + n q) m_1 , \notag
\end{align}
with $\lambda = 1-p-q$ as defined before.}
We can now compute $t_1$: after a few standard steps we obtain
\begin{align*}
    t_1 &= {\mathbb V}{\rm ar}\big[ A_n(1)A_n(2) \big]= \lambda^2 m_4 + \lambda \big(p-q+2nq\big) m_3+\big(nq(1-q)+n^2q^2\big)m_2 - \big(\mathfrak{m}_1\big)^2. 
\end{align*}

{\it Step 3: evaluation of the $t_k$ for $k=2,3,\ldots$} The rest of the proof focuses on the evaluation of $t_2,t_3,\ldots$\MM{;} together with our expression for $t_1$, this leads to an expression for $v_1$. To this end, we rewrite, for $k=2,3,\ldots$,
\[t_k= {\mathbb E}\Big[A_n(1)A_n(2)\Big({\mathbb E}\big[A_n(k)A_n(k+1)\,|\,A_n(2)\big]-{\mathbb E}\big[A_n(k)A_n(k+1)\big]\Big)\Big],\]
where we use that $\{A_n(k)\}_{k\in{\mathbb N}}$ is a Markov chain. 
Now notice that, conditional on $A_n(k)$, 
 \begin{align*}
  A_n(k+1) &\sim {\mathbb B}{\rm in}(A_n(k), r_1) + {\mathbb B}{\rm in}(n-A_n(k),s_1)\\&={\mathbb B}{\rm in}(A_n(k), 1-p) + {\mathbb B}{\rm in}(n-A_n(k),q),
 \end{align*}
 with the two binomial random variables appearing in the right-hand side being independent. 
This implies that
\begin{align*}
 {\mathbb E}\big[A_n(k)A_n(k+1)\,|\,A_n(2)\big]&=\lambda\,{\mathbb E}\big[(A_n(k))^2\,|\,A_n(2)]+ n q   \, {\mathbb E}\big[A_n(k)\,|\,A_n(2)\big].
\end{align*}
The next step is to use standard formulas for the first and second moment of the binomial variable, in combination with the distributional equality (for given $A_n(2)$)
 \[A_n(k) \sim {\mathbb B}{\rm in}(A_n(2), r_{k-1}) + {\mathbb B}{\rm in}(n-A_n(2),s_{k-1}),\]
 again with the binomial random variables in the right-hand side independent. 
We thus arrive, after a number of standard computations, at
\begin{equation}\label{aatjes} {\mathbb E}\big[A_n(k)A_n(k+1)\,|\,A_n(2)\big] =
\mathfrak{a}_{0,k}+ \mathfrak{a}_{1,k} \,A_n(2) + \mathfrak{a}_{2,k}\, (A_n(2))^2, \end{equation}
where
\begin{align*}
  \mathfrak{a}_{0,k} &:= \frac{n\varrho^2(n - 1)}{\lambda^3} \lambda^{2k} - \frac{n\varrho(\lambda - 2\lambda\varrho + n\varrho(1+\lambda))}{\lambda^2} \lambda^k + \mathfrak{m}_{1},\\
  \mathfrak{a}_{1,k} &:= \left(\frac{2\varrho - 1-2n\varrho}{\lambda^3}\right)\lambda^{2k} + \left(\frac{\lambda(1-  2\varrho) + n \varrho \, (1+\lambda)}{\lambda^2} \right)\lambda^k\\
  \mathfrak{a}_{2,k} &:= \left(\frac{1}{\lambda^3}\right)\lambda^{2k}.
\end{align*}

{\it Step 4: derivation of final expression.} 
Combining \eqref{aatjes} with the definition of the objects $\mathfrak m_i$ in \eqref{eq:m1}, \eqref{eq:m2}, and \eqref{eq:m3},
it can be verified that
\begin{align*}t_k&= {\mathbb E}\Big[A_n(1)A_n(2)\Big({\mathbb E}\big[A_n(k)A_n(k+1)\,|\,A_n(2)\big]-{\mathbb E}\big[A_n(k)A_n(k+1)\big]\Big)\Big]= \mathfrak{c}_1 \lambda^k + \mathfrak{c}_2 \lambda^{2k} \,,
\end{align*}
where
\begin{align*}
    & \mathfrak{c}_1 := \frac{\lambda(1-  2\varrho) + n \varrho \, (1+\lambda)}{\lambda^2} \, \mathfrak{m}_2 - \frac{n\MM{\lambda}\,\varrho (1-  2\varrho) + n^2\varrho^2\,(1+\lambda)}{\lambda^2} \, \mathfrak{m}_1,\\
    & \mathfrak{c}_2 := \frac{1}{\lambda^3} \, \mathfrak{m}_3 + \frac{2\varrho - 1-2n\varrho}{\lambda^3} \, \mathfrak{m}_2 +  \frac{n\varrho^2(n - 1)}{\lambda^3} \, \mathfrak{m}_1 .
\end{align*}
Upon combining the above steps, we thus find our final expression:
\begin{equation}\label{v1}
v_1 = t_1 + 2\sum_{k = 2}^{\infty} t_k = t_1 + 2 \left(\mathfrak{c}_1  \frac{\lambda^2}{1-\lambda} + \mathfrak{c}_2 \frac{\lambda^4}{1-\lambda^2} \right).
\end{equation}}

We finally compute $c_{01}$. \MRHM{Using standard computation rules for covariances of sums, in combination with the imposed stationarity, we need to evaluate
\begin{align}  c_{01}
=& \sum_{k = 1}^\infty {\mathbb C}{\rm ov}\big[A_n(1), \, A_n(k)A_n(k+1) \big]
 \label{ces}+
\sum_{k = 2}^\infty {\mathbb C}{\rm ov} \big[A_n(k), \, A_n(1)A_n(2) \big] 
. \end{align}
}
We first focus on the calculation of ${\mathbb C}{\rm ov}\big[A_n(1), \, A_n(k)A_n(k+1) \big]$,  now relying on  
\[A_n(k) \sim {\mathbb B}{\rm in}(A_n(1), r_{k}) + {\mathbb B}{\rm in}(n-A_n(1),s_{k}),\]
with the two quantities in the right-hand side independent.
It follows that
\[ {\mathbb E}\big[A_n(k)A_n(k+1)\,|\,A_n(1)\big]=\lambda\,{\mathbb E}\big[(A_n(k))^2\,|\,A_n(1)]+ n q   \, {\mathbb E}\big[A_n(k)\,|\,A_n(1)\big].\]
After some elementary \michel{calculations, using the same reasoning as in our derivation of $v_1$,} we conclude that
\begin{align*}
 {\mathbb E}\big[A_n(k)A_n(k+1)A_n(1)\big] &=  {\mathbb E}\big[A_n(1)\,{\mathbb E}\big[A_n(k)A_n(k+1)\,|\,A_n(1)\big] \big] = \mathfrak{d}_0 + \mathfrak{d}_1 \lambda^k + \mathfrak{d}_2 \lambda^{2k},
\end{align*}
where
\begin{align*}
    & \michel{\mathfrak{d}_0 := \mathfrak{m}_1 m_1} \\
    & \mathfrak{d}_1 := \frac{\lambda - 2\lambda\varrho + nq + 2\lambda n\varrho}{\lambda} \, m_2 - \frac{n\varrho(\lambda - 2\lambda\varrho + nq + 2\lambda n\varrho)}{\lambda} \, m_1 \\
    & \mathfrak{d}_2 := \frac{1}{\lambda} \, m_3 - \frac{2n\varrho - 2\varrho + 1}{\lambda} \, m_2 + \frac{n\varrho^2(n - 1)}{\lambda} \, m_1 \,\js{= 0} .
\end{align*}
\js{\MM{We conclude that} ${\mathbb C}{\rm ov}\big[A_n(1), \, A_n(k)A_n(k+1) \big]= \mathfrak{d}_1 \lambda^k $.}
To obtain ${\mathbb C}{\rm ov} \big[A_n(k), \, A_n(1)A_n(2) \big]$, we use the same recipe as the one developed when computing $v_1$. Indeed,
\begin{align*}
    {\mathbb C}{\rm ov} \big[A_n(k), \, A_n(1)A_n(2) \big] &= {\mathbb E}\big[A_n(1)A_n(2)A_n(k)\big] -{\mathbb E}\big[A_n(1)A_n(2)\big]{\mathbb E}\big[A_n(k)\big] \\
    & = {\mathbb E}\big[A_n(1)A_n(2) \, {\mathbb E}\big[A_n(k) \, | \,A_n(2) \big]\big] - {\mathbb E}\big[A_n(1)A_n(2)\big]{\mathbb E}\big[A_n(k)\big] \\
    & = \left(\mathfrak{m}_2 \, \frac{1}{\lambda^2}
    - \mathfrak{m}_1 
    \frac{n \varrho}{\lambda^2}\right) \lambda^k.
\end{align*}
\js{
Note that $\{A_n(k)\}_{k\in{\mathbb Z}}$ is, in stationarity, a {\it time-reversible} Markov chain, and as a consequence the vector $(A_n(1), A_n(2), A_n(k+1))$ has the same joint distribution as $(A_n(k+1), A_{n}(k), A_n(1))$.  \MM{This means ${\mathbb C}{\rm ov} \big[A_n(k+1), \, A_n(1)A_n(2) \big] = {\mathbb C}{\rm ov} \big[A_n(1), \, A_n(k)A_n(k+1) \big]$. In particular, it can be verified that $\mathfrak{d}_1 \lambda = \mathfrak{m}_2-\mathfrak{m}_1 n\varrho$.}
}
Upon combining the above, we finally obtain that the expression after the last equality sign in display \eqref{ces} equals
\begin{align}\notag
c_{01} &=
 \sum_{k = 1}^\infty {\mathbb C}{\rm ov}\big[A_n(1), \, A_n(k)A_n(k+1) \big]+\sum_{k = 2}^\infty {\mathbb C}{\rm ov} \big[A_n(k), \, A_n(1)A_n(2) \big] 
\\
& \js{ = \mathfrak{d}_1 \frac{\lambda}{1-\lambda} + \left(\mathfrak{m}_2 - \mathfrak{m}_1 n \varrho \right) \frac{1}{1-\lambda}} \:\MM{= 2\mathfrak{d}_1 \frac{\lambda}{1-\lambda}},\label{c01}
\end{align}

\js{
\begin{remark} \em In case $\lambda = 0$, i.e., $p+q=1$ (or $\varrho = q$), it was mentioned in Remark \ref{RPQ} that $\{A_n(k)\}_{k\in{\mathbb N}}$ is a sequence of {\it independent} random variables, with $A_n(k) \sim {\mathbb B}{\rm in}(n, \varrho)$. It follows that 
\[
\begin{aligned}
    & \nu_0 = n\varrho(1-\varrho), \\
    & \nu_1 = {\mathbb V}{\rm ar} \left[A_n(1)A_n(2)\right] = (n\varrho(1-\varrho)+n^2\varrho^2)^2 - (n^2\varrho^2)^2 ,\\
    & c_{01} = {\mathbb C}{\rm ov} \left[A_n(1), A_n(1)A_n(2)\right] = (n\varrho(1-\varrho)+n^2\varrho^2)(n\varrho) - (n\varrho)^3.
\end{aligned}
\]
These expressions are consistent with earlier findings. In particular, if $p+q=1$ then it follows from \eqref{eq:nu0} that $\nu_0 = n\varrho(1-\varrho)$. Also, \[\mathfrak{m}_1 = (n\varrho)^2 = (m_1)^2,\quad\quad\mathfrak{m}_2 = (n\varrho(1-\varrho)+n^2\varrho^2) \, n\varrho = m_1m_2.\] Substituting this into \eqref{v1} and \eqref{c01}, we have $\nu_1 = t_1 = (m_2)^2 - (m_1)^4$, and $c_{01} = \mathfrak{m}_2 - \mathfrak{m}_1 n \varrho = m_1m_2 - (m_1)^3$. This means that, despite the fact that we have been working under $\MM{\lambda} \ne 0$ in previous derivations, \eqref{eq:nu0}, \eqref{v1}, and \eqref{c01} hold for $\MM{\lambda}=0$. \hfill$\Diamond$ \end{remark}
}

\section{Verification of finite (co-)variances}\label{AppB}

 For the asymptotic normality to hold, one needs that $v_0$, $v_1$, and $c_{01}$ are finite. In this appendix, we point out how this can be proven for $v_0$ in the case that \michel{$X\sim{\mathbb P}{\rm ar}(C_X,\alpha)$ and $Y\sim{\mathbb P}{\rm ar}(C_Y,\beta)$, leading to the interpretable condition featuring in Proposition~\ref{P3}.} For $v_1$ and $c_{01}$ similar reasonings apply, but the underlying computations are considerably more cumbersome.
 Our argumentation essentially follows the argumentation of the proof of \cite[Theorem 1]{TAQ}, in combination with some elements from \cite{BM}. \michel{The proof does not cover the technically more demanding case of $\min\{\alpha,\beta\}$ being integer-valued.}

\medskip

{\it Step 1: evaluate the required probability generating function.}
 Let $r_k$ be defined as before, but with the on-time just having started at time $1$; likewise, $s_k$ is as before, but with the off-time just having started at time $1$. The following two equations are evident:
\begin{align*}
    r_k = \sum_{\ell=1}^{k-1} f_\ell \, s_{k-\ell}  + \sum_{\ell=k}^\infty f_\ell \,, \qquad s_k = \sum_{\ell=1}^{k-1} g_\ell \, r_{k-\ell}
\end{align*}
The next step is to find the probability generating functions pertaining to $r_k$ and $s_k$, using computations similar to those presented in \cite{BM}.
To this end, in both equations in the previous display, we multiply both sides by $z^k$, and define, for $|z|<1$,
\[F(z):= \sum_{k=1}^{\infty}z^k \, f_k,\:\:\:\:G(z):= \sum_{k=1}^{\infty}z^k \, g_k.\]  
We thus find, by swapping the summation order and splitting $z^k$ into $z^\ell z^{k-\ell}$ (or just by recognizing the convolution structure),
\begin{align*}
    S(z) &:= \sum_{k=1}^{\infty}z^k \, s_k= \sum_{k=1}^{\infty}z^k \sum_{\ell=1}^{k-1} g_\ell \, r_{k-\ell} = \sum_{\ell=1}^\infty z^\ell\sum_{k=\ell+1}^\infty z^{k-\ell} g_\ell r_{k-\ell}=G(z)\,R(z),
\end{align*}
where $R(z) := \sum_{k=1}^{\infty}z^k \, r_k.$
Likewise, 
\begin{align*}
R(z) &= F(z)\,S(z) + \sum_{k=1}^\infty z^k\sum_{\ell=k}^\infty f_\ell = F(z)\,S(z) + \sum_{\ell=1}^\infty  f_\ell  \sum_{k=1}^\ell z^k\\
    &=F(z)\,S(z) + z \sum_{\ell=1}^\infty  f_\ell  \frac{1-z^\ell}{1-z}=F(z)\,S(z) +H(z); \:\:\:\:\:\:\:H(z):= z     \frac{1-F(z)}{1-z}.
\end{align*}
After some elementary \michel{calculations}, we \michel{find}
\[
R(z) = \frac{H(z)}{1-F(z)\,G(z)}, \qquad S(z) = \frac{G(z)\,H(z)}{1-F(z)\,G(z)} \,.
\]
Let the probability $r_k^{\rm (res)}$ be as defined as $r_k$, but then with the residual on-time being distributed according to $(\bar f_\ell)_{\ell\in{\mathbb N}}$ rather than $(f_\ell)_{\ell\in{\mathbb N}}$. As a consequence, 
\begin{align*}
    r_k^{\rm (res)} = \sum_{\ell=1}^{k-1} \bar f_\ell \, s_{k-\ell}  + \sum_{\ell=k}^\infty \bar f_\ell, 
\end{align*}
and hence, in self-evident notation
\[R^{\rm (res)}(z)= \bar F(z)\,S(z) + z     \frac{1-\bar F(z)}{1-z}.\]
Also recall that $F(z)$ and $\bar F(z)$ can be expressed in one another:
\begin{align*}
\bar F(z) &= \frac{1}{{\mathbb E}\,X} \sum_{k=1}^\infty z^k \sum_{\ell=k}^\infty f_\ell = \frac{1}{{\mathbb E}\,X}\cdot z\sum_{\ell= 1}^\infty  \frac{1-z^\ell}{1-z} f_\ell = \frac{H(z)}{{\mathbb E}\,X}= \frac{H(z)}{F'(1)}.
\end{align*}
The next step is to combine the above, so as to determine, for $|z|<1$,
\begin{align}R_-^{\rm (res)}(z)&:= \sum_{k=1}^{\infty}z^k \, (r_k^{\rm (res)}-\varrho) = R^{\rm (res)}(z)-\varrho\frac{z}{1-z}\notag \\
    &= \frac{H(z)}{F'(1)}\left(\frac{G(z)\,H(z)}{1-F(z)\,G(z)}-\frac{z}{1-z}\right) +(1-\varrho)\frac{z}{1-z}\notag \\
    &=\frac{z}{1-z}\left(-\frac{(1-G(z))H(z)}{F'(1)(1-F(z)\,G(z))}+1-\varrho\right).\label{Rres}
\end{align}

\medskip

{\it Step 2: use a Tauberian theorem to characterize $r_k^{\rm (res)}-\varrho$ as $k\to\infty$.} 
\michel{The goal of this step is to use the expression we found for $R_-^{\rm (res)}(z)$ to show that $r_k^{\rm (res)}-\varrho$ vanishes as $k^{1-\alpha}$.
To do so, we expand \eqref{Rres} at $z\uparrow 1.$ From $X\sim{\mathbb P}{\rm ar}(C_X,\alpha)$, we have that
\[\lim_{i\to\infty} {\mathbb P}(X\geqslant i)\,i^\alpha = C_X^\alpha,\]
so that $ {\mathbb P}(X\geqslant i)$ is regularly varying with index $\alpha$ and the slowly varying function (as appearing in the definition of regular variation) equal to the constant $C_X^\alpha$. 
Tauberian theory \cite{BGT} entails that there is equivalence between (a)~the tail distribution of a non-negative random variable being regularly varying at $\infty$, and (b)~the associated probability generating function being regularly varying at $1$ (i.e., when approaching $1$ from below); see e.g.\ \cite[Theorem 8.1.6]{BGT}. This concretely means that if $X\sim{\mathbb P}{\rm ar}(C_X,\alpha)$, then the probability generating function of $X$ is of the form, as $z\downarrow 0$,
\[F(z) = \sum_{m=0}^n \gamma_m(1-z)^m + \gamma (1-z)^\alpha +o((1-z)^\alpha)\]
for some constants $\gamma_0,\ldots,\gamma_n$ and $\gamma$, with $n$ such that $\alpha\in(n,n+1)$; here $\gamma_0=1$ and $\gamma_1= -F'(1)=-{\mathbb E}\,X<0$.
In the sequel, we refer to this statement as `$F(z)$ is regularly varying of degree $\alpha$', or $F(z)\in {\mathbb R}{\rm V}_{\alpha}$ (where, in our case, the slowly varying function in the definition of regular variation is a constant). Observe that this definition can be used for any $F(z)$ that maps $|z|<1$ to real numbers, i.e., not just probability generating functions.  }

\michel{In this proof we consider the case that the tail of $Y$ is as most as heavy as the tail of $X$, in the sense that $\alpha\leqslant \beta$; the opposite case, in which $\beta\leqslant \alpha$, can be dealt with analogously. The fact that $Y\sim{\mathbb P}{\rm ar}(C_Y,\beta)$ entails that
\[\lim_{i\to\infty} {\mathbb P}(Y\geqslant i)\,i^\beta = C_Y^\beta.\]
 Applying the same argumentation as above, we can write, as $z\downarrow 0$,
\[G(z) = \sum_{m=0}^{n'} \delta_m(1-z)^m + \delta (1-z)^\beta+o((1-z)^\beta)\]
with $\beta\geqslant \alpha$ such that $\beta\in(n',n'+1)$ with $n'\geqslant n$, 
for some constants $\delta_0,\ldots,\delta_{n'}$ and $\delta$, 
or, in compact notation, $G(z)\in {\mathbb R}{\rm V}_{\beta}$; here $\delta_0=1$ and $\delta_1= -G'(1)=-{\mathbb E}\,Y<0$.}

\michel{Due to the fact that $\alpha\leqslant \beta$, we have, for some constants $\eta_1,\ldots\eta_n$ and $\eta$,
\begin{align*}
    \frac{1-F(z)\,G(z)}{1-z} &=\frac{1}{1-z}\Bigg(1-\left(1+\sum_{m=1}^n \gamma_m(1-z)^m + \gamma (1-z)^\alpha +o((1-z)^\alpha)\right)\\
    &\hspace{2.25cm}\left(1+\sum_{m=1}^{n'} \delta_m(1-z)^m + \delta (1-z)^\beta+o((1-z)^\beta)\right)\Bigg)\\
    &=\frac{1}{1-z}\left(\sum_{m=1}^{n}\eta_{m}(1-z)^m+\eta(1-z)^{\alpha}+o((1-z)^{\alpha})\right)\\
    &=\sum_{m=0}^{n-1}\eta_{m+1}(1-z)^m+\eta(1-z)^{\alpha-1}+o((1-z)^{\alpha-1})\in {\mathbb R}{\rm V}_{\alpha-1}
\end{align*}
as $z\uparrow 1$.
It is straightforward to see that $(1-G(z))/(1-z) \in  {\mathbb R}{\rm V}_{\beta-1}$. Again using $\alpha\leqslant \beta$, in combination with the definition of $H(z)$, we conclude that, for constants $\bar\eta_1,\ldots,\bar\eta_n$ and $\bar\eta$,
\begin{align*}
    I(z)&:=\frac{(1-G(z))H(z)}{F'(1)(1-F(z)\,G(z))}  =\frac{1}{F'(1)}\frac{1-G(z)}{1-z}\frac{1-z}{1-F(z)\,G(z)}H(z)\\
    &=\:\frac{1}{F'(1)}\left(-\sum_{m=0}^{n'-1} \delta_{m+1}(1-z)^m + \delta (1-z)^{\beta-1}\right)\left(\sum_{m=0}^{n-1}\bar\eta_{m+1}(1-z)^m+\bar\eta(1-z)^{\alpha-1}\right)\\
    &\hspace{2cm} z\left(\sum_{m=0}^n \gamma_{m+1}(1-z)^m +\gamma(1-z)^{\alpha-1}\right) + o((1-z))^{\alpha-1} \in {\mathbb R}{\rm V}_{\alpha-1}
\end{align*}
as $z\uparrow 1$.
As 
\[\lim_{z\uparrow 1}I(z)=\lim_{z\uparrow 1}\frac{(1-G(z))H(z)}{F'(1)(1-F(z)\,G(z))}=1-\varrho,\]
the expression between the brackets in \eqref{Rres} has no constant term. Hence all terms between the brackets in \eqref{Rres} are proportional to $(1-z)^m$, with $m\in\{1,\ldots,n-1\}$, or to $(1-z)^{\alpha-1}$, entailing that, as $z\uparrow 1$,
\[R_-^{\rm (res)}(z)=\frac{z}{1-z}\big(-I(z)+1-\varrho\,\big)\in {\mathbb R}{\rm V}_{\alpha-2}.\]
Now we can apply Tauberian theory in the reverse direction; to this end, we use that it can be argued that $r^{\rm (res)}_k-\varrho\geqslant 0$ for sufficiently large $k$, precisely as in the proof of \cite[Theorem 1]{TAQ}. We thus conclude that $r_k^{\rm (res)}-\varrho$ behaves as $k^{1-\alpha}$: for some positive constant $C$,
\[\left(r^{\rm (res)}_k-\varrho\right) \,k^{\alpha-1}\to C\]
as $k\to\infty.$}

\medskip

{\it Step 3: finding the range of $\alpha$ for which $v_0<\infty$.} In line with the calculations performed when computing $v_0$ for geometric on- and off times, 
\begin{align*}
        {\mathbb V}{\rm ar}\left[\sum_{k=1}^K A_n(k)\right] &= \sum_{k=1}^K{\mathbb V}\mbox{ar}\,A_n(k) + 2 \sum_{k=1}^K\sum_{k'=1}^{k-1} {\mathbb C}\mbox{ov}\left(A_n(k'),A_n(k)\right) \\
        &=  K\,{\mathbb V}\mbox{ar}\,A_n(1) + 2n \sum_{k=2}^{K}(K-k+1)\, {\mathbb C}{\rm ov}({\bs 1}(1),{\bs 1}(k)).
    \end{align*}
This means that
we have to verify for which $\alpha$ we have that 
\[\sum_{k=1}^{\infty}{\mathbb C}{\rm ov}({\bs 1}(1),{\bs 1}(k))=\varrho\sum_{k=1}^{\infty} \left(r^{\rm (res)}_k-\varrho\right)<\infty.\]
Using what was found in Step 2, we conclude that we need that $\alpha>2$. \michel{The analogous reasoning for $\beta\leqslant \alpha$ yields the condition $\beta>2$, so that we conclude that the claim holds if $\min\{\alpha,\beta\}>2.$}

\medskip

\end{document}